\documentclass[leqno]{amsart}

\usepackage{amsfonts}
\usepackage{amssymb, latexsym, amsmath, pb-diagram}
\usepackage{pstricks-add}
\usepackage{lamsarrow, pb-lams}
\usepackage{multicol}
\usepackage{bm}
\usepackage[mathscr]{eucal}
\usepackage{xy}
\usepackage{epic,eepic}
\xyoption{all}

\numberwithin{equation}{section}

\def\w{\widetilde}
\def\wh{\widehat}

\def\cf{{\it cf.}\ }
\def\o{\overline}

\newcommand{\RR}{\mathbb{R}}

\newcommand{\kk}{\Bbbk}
\newcommand{\ba}{\mathbf{a}}
\newcommand{\xx}{\mathbf{x}}
\newcommand{\yy}{\mathbf{y}}

\newcommand{\ZZ}{\mathcal{Z}}

\begin{document}
\title[The cohomology of generalized moment-angle complexes]
{\bf The homology coalgebra and cohomology algebra of generalized moment-angle complexes}
\author[Q. Zheng ]{Qibing Zheng}
\thanks{The author is supported by NSFC grant No. 10771105 and No. 11071125}
\keywords{generalized moment-angle complex, homotopy open cover, spectral sequence.}
\subjclass[2000]{Primary 13F55, 18G15, Secondary 16E05, 55U10.}
\address{School of Mathematical Science and LPMC, Nankai University,
Tianjin, 300071, P.R.China}
\email{zhengqb@nankai.edu.cn}
\maketitle

\begin{abstract}
In this paper, we compute the homology coalgebra and cohomology algebra over a field of all generalized moment-angle complexes and give a duality theorem on complementary moment-angle complexes.
\end{abstract}

The moment-angle complexes have been studied by topologists for many years
(\cf \cite{Porter} \cite{Lopez}). In 1990's
Davis and Januszkiewicz \cite{DJ} introduced toric manifolds which are studied
intensively by algebraic geometers. They observed
that every quasi-toric manifold is the quotient of a moment-angle complex $\mathcal{Z}_K$ by the free action
of a real torus. The topology of $\ZZ_K$ is complicated and getting more attentions by topologists lately
(\cf \cite{GM} \cite{Hoc} \cite{Baskakov} \cite{Panov} \cite{Franz}).
Recently a lot of work has been done on generalizing the moment-angle complex
$\ZZ_K=\ZZ_K(D^2,S^1)$ to space pairs $(X, A)$
(\cf  \cite{BB}, \cite{BBCG08}, \cite{GTh}, \cite{LP}). But in general, the (co)homology of a generalized moment-angle complex is not known. \vspace{3mm}

In this paper, we construct a (co)chain complex to compute the (co)homology group of a generalized moment-angle complex in Theorem 1.3 and give a duality theorem for a special type of generalized moment-angle complexes in Theorem 1.6. To compute the coalgebra structure of the homology and the algebra structure of cohomology, we construct two spectral sequences in Theorem 2.3 that converge respectively to the homology coalgebra and cohomology algebra of a generalized moment-angle complex. By using the spectral sequence that collapse for all generalized moment-angle complexes, we determine the homology coalgebra and cohomology algebra structure of a generalized moment-angle complex in Theorem 3.8. We end the paper with an example.\vspace{3mm}

\section{The (co)homology group of generalized moment-angle complexes}\vspace{3mm}

{\bf Notations and Conventions}\, In this paper, $\kk$ is a field. All objects (groups, (co)chain complexes, (co)algebras, etc.) are vector spaces over $\kk$ and $\otimes$ means $\otimes_\kk$. So base and dual objects always exist and $H^*(-)={\rm Hom}_\kk(H_*(-),\kk)$ and $H_*(C_1{\otimes}C_2,d_1{\otimes}d_2)=H_*(C_1,d_1){\otimes}H_*(C_2,d_2)$. For a topological space or a simplicial complex $X$,  $H_*(X)=H_*(X;\kk)$ and $H^*(X)=H^*(X;\kk)$, the (co)homology over $\kk$.

For a positive integer $n$, $[n]$ denotes the set $\{1,2,\cdots,n\}$.

A simplicial complex $K$ with vertex set $S$ is a set of subsets of $S$ such that if $\tau\in K$ and $\sigma\subset\tau$, then $\sigma\in K$ and $\{s\}\in K$ for every $s\in S$. $2^S$ is the full simplicial complex consisting of all subsets of $S$. The geometrical realization of $K$ is denoted by $|K|$. Notice that all simplicial complexes have the empty set $\phi$ as a $-1$-dimensional simplex except the void simplicial complex $\{\}$ that has no simplex at all.\vspace{3mm}

{\bf Definition 1.1} {\it Let $(X,A)=\{(X_k,A_k)\}_{k=1}^m$ be a sequence of topological space pair. For a simplicial complex $K$ with vertex set a subset of $[m]$, the topological space $\ZZ_K(X,A)$ is defined as follows. For a subset $\sigma$ of $[m]$, define
\[D(\sigma)=Y_1{\times}\cdots{\times}Y_n,\quad Y_i=\left\{\begin{array}{cl}
X_i&{\rm if}\,\,i\in \sigma, \\
A_i&{\rm if}\,\,i\not\in \sigma.
\end{array}
\right.\]
Then $\ZZ_K(X,A)=\cup_{\sigma\in K}D(\sigma)$.

In this paper, a generalized moment-angle complex $M=\ZZ_K(X,A)$ is a topological space defined as above that  satisfies the following condition. Every $A_k$ is either an open subspace of $X_k$ or a deformation retract of an open subspace of $X_k$. Specifically, if every $(X_k,A_k)$ is a CW-complex pair, $\ZZ_K(X,A)$ is a generalized moment-angle complex.}\vspace{3mm}

Notice that we do not require that the vertex set of $K$ is $[m]$. This will give a complete conclusion for duality.\vspace{3mm}

{\bf Definition 1.2} {\it Let $K$ be a simplicial complex. For a simplex $\sigma\in K$, the link and star of $\sigma$ are respectively the simplicial complex ${\rm link}_K(\sigma)=\{\eta\in K\,|\,\eta{\cup}\sigma\in K,\,\eta{\cap}\sigma=\phi\}$ and ${\rm star}_K(\sigma)=\{\tau\in K\,|\,\sigma{\cup}\tau\in K\}$. For a subset $\omega$ of the vertex set of $K$, $K|_\omega=\{\eta{\cap}\omega\,|\,\eta\in K\}$. For $\sigma\in K$ and $\sigma{\cap}\omega=\phi$, $K_{\sigma,\omega}={\rm link}_K(\sigma)|_\omega=\{\omega{\cap}\eta\,|\,\eta\in {\rm link}_K(\sigma)\}$.

For a generalized moment-angle complex $M=\ZZ_K(X,A)$, $(X,A)=\{(X_k,A_k)\}_{k=1}^m$, let $i_k\colon H_*(A_k)\to H_*(X_k)$ be the homomorphism induced by the inclusion map and $\frak k_k={\rm ker}i_k$, $\frak i_k={\rm coim}i_k={\rm im}i_k$, $\frak c_k={\rm coker}i_k$. $\Sigma=\{k\,|\,\frak c_k\neq 0\}$, $\Omega=\{k\,|\,\frak k_k\neq 0\}$, $I_M=\{(\sigma,\omega)\,|\,\sigma\in K,\,\sigma\subset\Sigma,\,\omega\subset\Omega,\,\sigma{\cap}\omega=\phi\}$. For $(\sigma,\omega)\in I_M$, $T_*^{\sigma,\omega}(M)=\frak a_1{\otimes}{\cdots}{\otimes}\frak a_m$, where $\frak a_k=\frak c_k$ if $k\in\sigma$, $\frak a_k=\frak k_k$ if $k\in\omega$, $\frak a_k=\frak i_k$ otherwise.
Dually, let $(-)^*$ denote the dual map and space and so $\frak k^*_k={\rm coker}i^*_k$, $\frak i^*_k={\rm im}i^*_k={\rm coim}i^*_k$, $\frak c^*_k={\rm ker}i^*_k$. For $(\sigma,\omega)\in I_M$, let $T^*_{\sigma,\omega}(M)=\frak a_1{\otimes}{\cdots}{\otimes}\frak a_m$, where $\frak a_k=\frak c^*_k$ if $k\in\sigma$, $\frak a_k=\frak k^*_k$ if $k\in\omega$, $\frak a_k=\frak i^*_k$ otherwise. The relation of the above graded groups are as in the  following diagram
\[\begin{array}{ccc}
H_*(A_k) &\stackrel{i_k}{\longrightarrow}&H_*(X_k)\\
\|&&\|\\
\frak k_k{\oplus}\frak i_k&\longrightarrow&\frak i_k{\oplus}\frak c_k
  \end{array}\quad\quad
  \begin{array}{ccc}
H^*(X_k) &\stackrel{i^*_k}{\longrightarrow}&H^*(A_k)\\
\|&&\|\\
\frak c^*_k{\oplus}\frak i^*_k&\longrightarrow&\frak i^*_k{\oplus}\frak k^*_k.
  \end{array}\]
}\vspace{3mm}

{\bf Theorem 1.3} {\it Let  $M=\ZZ_K(X,A)$, $(X,A)=\{(X_k,A_k)\}_{k=1}^m$, be a generalized moment-angle complex and everything is as in Definition 1.2. Then there are group isomorphisms
\[H_{*}(M)\cong\oplus_{(\sigma,\omega)\in I_M}H_*^{\sigma,\omega}(M),\,H_{k+1}^{\sigma,\omega}(M)=\oplus_{s+t=k}\w H_{s}(K_{\sigma,\omega}){\otimes}T_t^{\sigma,\omega}(M),\]
\[H^{*}(M)\cong\oplus_{(\sigma,\omega)\in I_M}H^*_{\sigma,\omega}(M)\,,\,H^{k+1}_{\sigma,\omega}(M)=\oplus_{s+t=k}\w H^{s}(K_{\sigma,\omega}){\otimes}T^t_{\sigma,\omega}(M).\]}

\begin{proof} We may suppose that every $A_k$ is an open subspace of $X_k$, otherwise, replace $A_k$ by the open subspace of $X_k$ of which $A_k$ is a deformation retract.

Construct chain complex $(C_*(M),d)$ as follows. Let $(U_k,d_k)$ be the chain complex such that $U_k=\frak k_k{\oplus}\frak i_k{\oplus}\frak c_k\oplus\frak q_k$, $d_k(x)=0$ for all $x\in\frak k_k{\oplus}\frak i_k{\oplus}\frak c_k$  and $d_k\colon\frak q_k\to \frak k_k$ is an isomorphism with degree lowered by $1$ if $\frak k_k\neq 0$ and $\frak q_k=0$ if $\frak k_k=0$. Let $(S_k,d_k)$ be the chain subcomplex of $(U_k,d_k)$ such that $S_k=\frak k_k{\oplus}\frak i_k$. For a simplex $\sigma\in K$, $(C_*(\sigma),d)=(V_1{\otimes}{\cdots}{\otimes}V_m,d_1{\otimes}{\cdots}{\otimes}d_m)$, where $V_k=U_k$ if $k\in\sigma$ and $V_k=S_k$ if $k\not\in\sigma$. Then for $\sigma\subset\sigma'$, $(C_*(\sigma),d)$ is a chain subcomplex of $(C_*(\sigma'),d)$ and for all $\sigma,\tau\in K$, $C_*(\sigma){\cap}C_*(\tau)=C_*(\sigma{\cap}\tau)$. Define $(C_*(M),d)$ to be the chain complex such that $C_*(M)=+_{\sigma\in K}C_*(\sigma)$ (not direct sum $\oplus$!).

Let $(S_*(X),d)$ be the singular chain complex over $\kk$ of the topological space $X$. Take a fixed representative in $S_*(A_k)$ for all homology classes in a base of $\frak k_k,\frak i_k$, a fixed representative in $S_*(X_k)$ for all homology classes in a base of $\frak c_k$ and a base in $S_*(X_k)$ for $\frak q_k$ such that for every base element $y\in\frak q_k$, $dy$ is the representative for the base element $d_ky$. Then we may regard $(S_k,d_k)$ and $(U_k,d_k)$ respectively as chain subcomplex of $(S_*(A_k),d)$ and $(S_*(X_k),d)$. The two inclusion homomorphisms of chain complexes are homotopy equivalences and so induce homology isomorphism.

Now we prove $H_*(C_*(M),d)=H_*(M)$. We use double induction on the dimension and the number of maximal simplices (a maximal simplex is not the proper face of any simplex) of $K$. If dim$K=-1$, the conclusion is trivial. If $K$ has only one maximal simplex $\sigma$, i.e., $K=2^\sigma$, then there is a chain complex homotopy equivalence \[i_M\colon(C_*(M),d)=(C_*(\sigma),d)\simeq(S_*(Y_1){\otimes}{\cdots}{\otimes}S_*(Y_m),d{\otimes}{\cdots}{\otimes}d)
\simeq (S_*(Y_1{\times}{\cdots}{\times}Y_m),d)=(S_*(D(\sigma)),d),\]
where $D(\sigma)$ is as defined in Definition 1.1, the first chain homotopy equivalence is the inclusion homomorphism and the second is the Eilenberg-Zilber chain homotopy equivalence. Suppose there is a chain homotopy equivalence $i_M\colon (C_*(M),d)\to (S_*(M),d)$ for all moment-angle complex $M=\ZZ_K(X,A)$ such that $K$ has $\leqslant u$ maximal simplices and that dim$K\leqslant v$. Then for $K$ with maximal simplices $\sigma_1,\cdots,\sigma_{u+1}$ and dim$K\leqslant v$, denote by $K_1$ the simplicial complex with maximal simplices $\sigma_1,\cdots,\sigma_u$, $K_2=2^{\sigma_{u+1}}$ and $K_3=K_1{\cap}K_2$. Denote $M_i=\ZZ_{K_i}(X,A)$. Then by induction hypothesis, there is a homotopy equivalence $i_{M_k}\colon (C_*(M_k),d)\to (S_*(M_k),d)$ for $k=1,2,3$.
From the following commutative diagram
\[\begin{array}{ccccccc}
0\to&C_*(M_3)&\to&C_*(M_1){\oplus}C_*(M_2)&\to&C_*(M)&\to 0\\
&{\scriptstyle i_{M_3}}\downarrow\,\,\quad&&{\scriptstyle i_{M_1}\oplus i_{M_2}}\downarrow\quad\quad\,\,\,&&{\scriptstyle i'_{M}}\downarrow\,\,\quad\\
0\to&S_*(M_3)&\to&S_*(M_1){\oplus}S_*(M_2)&\to&(S_*(M_1){\oplus}S_*(M_2))/S_*(M_3)&\to 0
  \end{array}\]
we get a homotopy equivalence $i'_M$. By excision axiom, there is a homotopy equivalence $j\colon(S_*(M),d)\to ((S_*(M_1){\oplus}S_*(M_2))/S_*(M_3),\w d)$ and so $i_M=j^{-1}i'_M$ is a homotopy equivalence. For $K$ with maximal simplices $\sigma_1,\cdots,\sigma_{u}$ and dim$K=v{+}1$, we may suppose $|\sigma_u|=v{+}1$. Denote by $K_1$ the simplicial complex with maximal simplices $\sigma_1,\cdots,\sigma_{u-1}$, $K_2=2^{\sigma_{u}}$ and $K_3=K_1{\cap}K_2$. Then $K_1,K_2,K_3$ all satisfy the induction hypothesis and the homotopy equivalence $i_M$ exists by the same induction proof as above. Thus, $H_*(C_*(M),d)=H_*(M)$ for all $M$.

Let $(C_*(\w X),d)=(U_1{\otimes}{\cdots}{\otimes}U_m,d_1{\otimes}{\cdots}{\otimes}d_m)$ and denote the element of $C_*(\w X)$ by $\ba=(a_1,{\cdots},a_m)$ with every $a_k\in U_k$. For $(a_1,\cdots,a_m)\in C_*(\w X)$, define\\
\hspace*{31mm}$\rho(\ba)=\{k\,|\,k\in\frak c_k\}, \quad \rho^*(\ba)=\{k\,|\,k\in\frak i_k\},\quad \varrho(\ba)=\{k\,|\,a_k\in\frak k_k\,{\rm or}\,\frak q_k\}$\\
and $C_*^{\sigma,\w\sigma}(\w X)=\{\ba\in C_*(\w X)\,|\,\rho(\ba)=\sigma,\,\rho^*(\ba)=\w\sigma\}$.

For $(\sigma,\omega)\in I_M$, $C_*^{\sigma,\omega}(M)=\{\ba\in C_*(M)\,|\,\rho(\ba)=\sigma,\,\varrho(\ba)=\omega\}$ is a chain subcomplex of $(C_*(M),d)$ such that $(C_*(M),d)=\oplus_{(\sigma,\omega)\in I_M}(C_*^{\sigma,\omega}(M),d)$. Let $(\w C_*(K_{\sigma,\omega}),d)$ be the augmented simplicial chain complex over $\kk$ of $K_{\sigma,\omega}$. Then there is a chain complex isomorphism $\psi\colon(C_{*+1}^{\sigma,\omega}(M),d)\to (\w C_*(K_{\sigma,\omega}),d){\otimes}T_*^{\sigma,\omega}(M)$ defined as follows. For $\ba=(a_1,\cdots,a_m)\in C_*^{\sigma,\omega}(M)$, let $\w\ba=(b_1,\cdots,b_m)$, where $b_k=a_k$ if $a_k\not\in\frak q_k$ and $b_k=d_k(a_k)$ if $a_k\in\frak q_k$. Then $\psi(\ba)=\eta{\otimes}\w\ba$ with $\eta=\{k\,|\,a_k\in\frak q_k\}$. So $H_{*+1}(C_*^{\sigma,\omega}(M),d)=\w H_*(K_{\sigma,\omega}){\otimes}T_*^{\sigma,\omega}(M)$.

Dually, we have\\
$\hspace*{5mm}H^*(M)={\rm Hom}_\kk(\oplus_{(\sigma,\omega)\in I_M}H_*^{\sigma,\omega}(M),\kk)
=\oplus_{(\sigma,\omega)\in I_M}{\rm Hom}_\kk(H_*^{\sigma,\omega}(M),\kk)
=\oplus_{(\sigma,\omega)\in I_M}H^*_{\sigma,\omega}(M)$,

$\hspace*{-2mm}H^{k+1}_{\sigma,\omega}(M)={\rm Hom}_\kk(H_{k+1}^{\sigma,\omega}(M),\kk)
={\rm Hom}_\kk(\oplus_{s+t=k}H_s(K_{\sigma,\omega}){\otimes}T_t^{\sigma,\omega}(M),\kk)
=\oplus_{s+t=k}H^s(K_{\sigma,\omega}){\otimes}T^t_{\sigma,\omega}(M)$.
\end{proof}

{\bf Definition 1.4} {\it Let $K$ be a simplicial complex with vertex set a subset of $[m]$. The Alexander dual of $K$ relative to $[m]$ is the simplicial complex $K^*=\{\o\sigma\,|\,\sigma\in K^c\}$, where $K^c=2^{[m]}{\setminus}K$ and $\o\sigma=[m]{\setminus}\sigma$. It is obvious that $(K^*)^*=K$.}\vspace{3mm}

By the above definition, the Alexander dual of $\{\phi\}$ relative to $[m]$ is $2^{[m]}{\setminus}\{[m]\}$, the Alexander dual of $\{\}$ relative to $[m]$ is $2^{[m]}$. $\ZZ_{\{\}}(X,A)=\phi$ and $\ZZ_{\{\phi\}}(X,A)=A_1{\times}{\cdots}{\times}A_m$.\vspace{3mm}

{\bf Theorem 1.5} {\it Let $M=\ZZ_K(X,A)$ be the topological space in Definition 1.1 (may not be a generalized moment-angle complex). Denote the topological space $M^c=\ZZ_{K^*}(X,B)$, $(X,B)=\{(X_k,B_k)\}_{k=1}^m$, where $B_k=X_k{\setminus}A_k$ and $K^*$ is as in Definition 1.4. Let $\w X=X_1{\times}{\cdots}{\times}X_m$. Then $M^c=\w X{\setminus}M$.}

\begin{proof} For a subset $\sigma$ of $[m]$, $2^\sigma$ has one maximal simplex $\sigma$. If $\sigma=\{i_0,\cdots,i_s\}$, then $(2^\sigma)^*$ has $s{+}1$ maximal simplices $\o{i_0},\cdots,\o{i_s}$, where $\o i=[m]{\setminus}\{i\}$. So $\w X{\setminus}D(\sigma)=\cup_{i\in\sigma}D^*(\o i)$, where $D^*$ is defined by that
\[D^*(\omega)=Z_1{\times}\cdots{\times}Z_n,\quad Z_i=\left\{\begin{array}{cl}
X_i&{\rm if}\,\,i\in \omega, \\
B_i&{\rm if}\,\,i\not\in \omega.
\end{array}
\right.\]
This implies that when $K=2^\sigma$, $\w X{\setminus}\ZZ_K(X,A)=\ZZ_{K^*}(X,B)$. It is easy to check that for any two simplicial complexes $K$ and $L$, $(K\cup L)^*=K^*\cap L^*$ and $\ZZ_{K^*\cap L^*}(X,B)=\ZZ_{K^*}(X,B)\cap\ZZ_{L^*}(X,B)$. So for any simplicial complex $K$,
\begin{eqnarray*}&&\quad\ZZ_{K^*}(X,B)\\
&&=\ZZ_{(\cup_{\sigma\in K}2^\sigma)^*}(X,B)=\ZZ_{\cap_{\sigma\in K}(2^\sigma)^*}(X,B)=\cap_{\sigma\in K}\ZZ_{(2^\sigma)^*}(X,B)\\
&&=\cap_{\sigma\in K}\big(\w X{\setminus}\ZZ_{2^\sigma}(X,A)\big)=\cap_{\sigma\in K}\big(\w X{\setminus}D(\sigma)\big)
=\w X{\setminus}\big(\cup_{\sigma\in K}D(\sigma)\big)\quad\\
&&=\w X{\setminus}\ZZ_K(X,A).
\end{eqnarray*}
\end{proof}

{\bf Theorem 1.6} {\it Let $M=\ZZ_K(X,A)$, $(X,A)=\{(X_k,A_k)\}_{k=1}^m$ be a generalized moment-angle complex satisfying the following conditions.

1) Every $X_k$ is a closed orientable manifold (with respect to homology over $\kk$) of dimension $r_k$.

2) Every $A_k$ is a polyhedron subspace of $X_k$ that is the deformation retract of a neighborhood.

Then $M^c=\ZZ_{K^*}(X,B)$ as in Definition 1.5 is also a generalized moment-angle complex and for any $(\sigma,\omega)\in I_M$ such that $\omega\neq\phi$, there are group  isomorphisms
\[\zeta_{\sigma,\omega}\colon H_*^{\sigma,\omega}(M)\cong H^{r-*-1}_{\w\sigma,\omega}(M^c),\quad
\zeta^*_{\sigma,\omega}\colon H^*_{\sigma,\omega}(M)\cong H_{r-*-1}^{\w\sigma,\omega}(M^c),\]
where $r=\Sigma_{k=1}^m r_k$, $\w\sigma=[m]\setminus(\sigma{\cup}\omega)$, $\zeta^*_{\sigma,\omega}$ is the dual map of $\zeta_{\sigma,\omega}$.

Precisely, the isomorphism $\zeta_{\sigma,\omega}$ is defined as follows. For $\omega\neq\phi$, we have an isomorphism $\o \zeta_{\sigma,\omega}$ that satisfies the following commutative diagram
\[\begin{array}{ccc}\w H_{s-1}(K_{\sigma,\omega})&
\stackrel{\partial^{-1}_{\sigma,\omega}}{-\!\!\!-\!\!\!\longrightarrow}&
H_s(2^{\omega},K_{\sigma,\omega})\\
&&\\
&_{\o\zeta_{\sigma,\omega}}\!\!\!\searrow\quad\quad&{\scriptstyle c_{\sigma,\omega}}\downarrow\quad\,\,\\
&&\\
&&\w H^{t-s}(K^*_{\w\sigma,\omega})
\end{array}\]
where $t$, $c_{\sigma,\omega}$ and $\partial_{\sigma,\omega}$ are defined in the proof of the theorem. For $x_k\in\frak k_k$, denote $\w x_k=\gamma_k(\partial^{-1}_k(x_k))$, where $\partial_k\colon H_{*+1}(X_k,A_k)\to H_*(A_k)$ is the connecting homomorphism and $\gamma_k\colon H_*(X_k,A_k)\to H^{r_k-*}(B_k)$ is the Alexander duality isomorphism. For $x_k\in\frak i_k$ or $\frak c_k$, denote $\w x_k=\gamma_k(x_k)$ with $\gamma_k\colon H_*(X_k)\to H^{r_k-*}(X_k)$ the Poncar\'{e} duality isomorphism. Then for
$a{\otimes}(x_1,\cdots,x_m)\in H_*^{\sigma,\omega}(M)$ with $a\in H_*(K_{\sigma,\omega})$ and $(x_1,\cdots,x_m)\in T_*^{\sigma,\omega}(M)$, $\zeta_{\sigma,\omega}\big(a{\otimes}(x_1,\cdots,x_m)\big)
=\o\zeta_{\sigma,\omega}(a){\otimes}(\w x_1,\cdots,\w x_m)$.
}
\begin{proof} Let everything be as in the proof of Theorem 1.3 and define
$(C_*(\w X,M),d)=(C_*(\w X)/C_*(M),\w d\,)$.

For $\eta\in K^*$, define $(\w C_*(\eta),d)=(W_1{\otimes}{\cdots}{\otimes}W_m,d_1{\otimes}{\cdots}{\otimes}d_m)$, where $(W_i,d_i)=(U_i,d_i)$ if $i\in \eta$ and
$(W_i,d_i)=(U_i/S_i,\w d_k\,)$ if $i\not\in \eta$. We prove that  $(C_*(\w X,M),d)=+_{\eta\in K^*}(\w C_*(\eta),d)$.

For  $\sigma=\{i_0,\cdots,i_s\}$, $(2^\sigma)^*$ has $s{+}1$ maximal simplices $\o{i_0},\cdots,\o{i_s}$ and so $C_*(\w X)/C_*(\sigma)=+_{i\in\sigma}\w C_*(\o i)$. This implies that $(C_*(\w X)/C_*(M_\sigma),\w d\,)=+_{\eta\in (2^\sigma)^*}(\w C_*(\eta),d)$, where $M_\sigma=\ZZ_{2^\sigma}(X,A)$. So for $K$,
\begin{eqnarray*}&&\quad+_{\eta\in K^*}(\w C_*(\eta),d)\\
&&=+_{\eta\in(\cup_{\sigma\in K}2^\sigma)^*}(\w C_*(\eta),d)=+_{\eta\in(\cap_{\sigma\in K}(2^\sigma)^*)}(\w C_*(\eta),d)=\cap_{\sigma\in K}\big(+_{\eta\in (2^\sigma)^*}(\w C_*(\eta),d)\big)\\
&&=\cap_{\sigma\in K}\big(C_*(\w X)/C_*(M_\sigma),\w d\,\big)=\cap_{\sigma\in K}\big(C_*(\w X)/C_*(\sigma),\w d\,\big)
=\big(C_*(\w X)/(+_{\sigma\in K}C_*(\sigma)),\w d\,\big)\quad\\
&&=\big(C_*(\w X)/C_*(M),\w d\,\big).
\end{eqnarray*}

So $(C_*(\w X,M),d)=+_{\eta\in K^*}(\w C_*(\eta),d)$. For $(\sigma,\omega)\in I_{M}$ and $\w\sigma=[m]{\setminus}(\sigma{\cup}\omega)$, define $C_*^{\w\sigma,\omega}(\w X,M)=\{\ba\in C_*(\w X,M)\,|\,\rho^*(\ba)\!=\!\w\sigma,\,\varrho(\ba)\!=\!\omega\}$ and $T_*^{\w\sigma,\omega}(\w X,M)=\{\ba\!=\!(a_1,\cdots,a_m)\in C_*^{\w\sigma,\omega}(\w X,M)\,|\,a_k\in\frak q_k$ for all $k\in\varrho(\ba)\}$. Then $(C_*(\w X,M),d)=\oplus_{(\sigma,\omega)\in I_{M}}(C_*^{\w\sigma,\omega}(\w X,M),d)$ and so for every $(\sigma,\omega)\in I_M$, $(C_*^{\sigma,\w\sigma}(\w X)/C_*^{\sigma,\omega}(M),d)=(C_*^{\w\sigma,\omega}(\w X,M),d)$.

Let $(K_{\sigma,\omega})^*$ be the Alexander dual of $K_{\sigma,\omega}$ relative to $\omega$ and $K^*$ be the Alexander dual of $K$ relative to $[m]$. For $\eta\in(K_{\sigma,\omega})^*$, i.e., $\eta\subset\omega$ and $\omega{\setminus}\eta\not\in K_{\sigma,\omega}$, i.e., $\eta\subset\omega$ and $\sigma{\cup}(\omega{\setminus}\eta)\not\in K$, we have $[m]\setminus(\sigma{\cup}(\omega{\setminus}\eta))=\w\sigma{\cup}\eta\in K^*$, i.e., $\eta\in(K^*)_{\w\sigma,\omega}$. So $(K_{\sigma,\omega})^*\subset(K^*)_{\w\sigma,\omega}$. Conversely, for $\eta\in(K^*)_{\w\sigma,\omega}$, i.e., $\eta\subset\omega$ and $\w\sigma{\cup}\eta\in K^*$, we have $[m]\setminus(\w\sigma{\cup}\eta)=\sigma{\cup}(\omega{\setminus}\eta)\not\in K$, i.e., $\omega{\setminus}\eta\not\in K_{\sigma,\omega}$, $\eta\in(K_{\sigma,\omega})^*$. So $(K^*)_{\w\sigma,\omega}\subset(K_{\sigma,\omega})^*$. Thus, $(K_{\sigma,\omega})^*=(K^*)_{\w\sigma,\omega}$ and we may denote $K^*_{\w\sigma,\omega}=(K_{\sigma,\omega})^*=(K^*)_{\w\sigma,\omega}$. The correspondence $\eta\to \omega{\setminus}\eta$ induces a dual complex isomorphism $(C_*(2^{\omega},K_{\sigma,\omega}),d)\cong (\w C^{t-*}(K^*_{\w\sigma,\omega}),\delta)$, where $t$ is the cardinality of $\omega$. Denote the corresponding induced homology isomorphism by $c_{\sigma,\omega}\colon H_*(2^{\omega},K_{\sigma,\omega})\cong \w H^{t-*}(K^*_{\w\sigma,\omega})$

There is a dual complex isomorphism $\w\psi\colon(C_{*+1}^{\w\sigma,\omega}(\w X,M),d)\to (\w C^*(K^*_{\w\sigma,\omega}),\delta){\oslash}T_*^{\w\sigma,\omega}(\w X,M)$ ($(A_*{\oslash}B_*)_k=\oplus_{t-s=k}A_s{\otimes}B_t$) defined as follows. For $\ba=(a_1,\cdots,a_m)\in C_*^{\w\sigma,\omega}(\w X,M)$, let $\w\ba=(b_1,\cdots,b_m)$, where $b_k=a_k$ if $a_k\not\in\frak k_k$ and $d_k(b_k)=a_k$ if $a_k\in\frak k_k$. Then $\w\psi(\ba)=\eta{\oslash}\w\ba$ with $\eta=\{k\,|\,a_k\in\frak k_k\}$. So $H_{*+1}^{\w\sigma,\omega}(\w X,M)=H_{*+1}(C_*^{\w\sigma,\omega}(\w X,M),d)= \w H^*(K^*_{\w\sigma,\omega}){\oslash}T_*^{\w\sigma,\omega}(\w X,M)$.

Thus, the long exact sequence
\[\begin{array}{ccccccccc}
{\cdots}\to & H_n(C_*(M),d)&\stackrel{i}{\longrightarrow}&H_n(C_*(\w X),d)&\stackrel{j}{\longrightarrow}&
H_n(C_*(\w X,M),\w d\,)&\stackrel{\partial}{\longrightarrow}&H_{n-1}(C_*(M),d)&\to {\cdots}
  \end{array}\]
satisfies that ker$i=\oplus_{(\sigma,\omega)\in I_M,\omega\neq\phi}H_*^{\sigma,\omega}(M)$, coker$j=\oplus_{(\sigma,\omega)\in I_{M},\omega\neq\phi}H_*^{\w\sigma,\omega}(\w X,M)$ and the long exact sequence is the direct sum of the following diagram for all $(\sigma,\omega)\in I_M$ and base element $\ba=(a_1,\cdots,a_m)\in T_*^{\sigma,\omega}(M)$,
\[\begin{array}{ccccccccc}{\cdots}\to&\w H_s(K_{\sigma,\omega}){\otimes}\ba&\to&
\w H_s(2^{\omega}){\otimes}\ba&\to&
H_s(2^{\omega},K_{\sigma,\omega}){\otimes}\ba&
\stackrel{\partial_{\sigma,\omega}}{-\!\!\!-\!\!\!\longrightarrow}&
\w H_{s-1}(K_{\sigma,\omega}){\otimes}\ba&\to {\cdots}\\
&&&&&&&&\\
&&&&&c_{\sigma,\omega}{\|}{\wr}&^{\partial'_{\sigma,\omega}}\!\!\!\nearrow&&\\
&&&&&&&&\\
&&&&&\w H^{t-s}(K^*_{\w\sigma,\omega}){\oslash}\ba'&&&
\end{array}\]
where $\ba'=(a'_1,\cdots,a'_m)$ with $d_k(a'_k)=a_k$ if $a_k\in\frak k_k$ and $a'_k=a_k$ otherwise. When $\omega\neq\phi$, $\w H_*(2^{\omega})=0$. We have an isomorphism $\partial'_{\sigma,\omega}\colon H_{*}^{\w\sigma,\omega}(\w X,M)\to H_{*-1}^{\sigma,\omega}(M)$.

The Alexander duality isomorphism $\gamma\colon H_*(\w X,M)\cong H^{r-*}(M^c)$ is also the direct sum of isomorphisms $\gamma_{\w\sigma,\omega}\colon H_*^{\w\sigma,\omega}(\w X,M)\cong H^{r-*}_{\w\sigma,\omega}(M^c)$ for all $(\sigma,\omega)\in I_{M}$. From the  following commutative diagram for $Y_k=B_k$ or $\phi$,
\[\begin{array}{ccc}
H_*(X_1,X_1{\setminus}Y_1){\otimes}{\cdots}{\otimes}H_*(X_m,X_m{\setminus}Y_m)&
\stackrel{\gamma_1{\otimes}{\cdots}{\otimes}\gamma_m}{-\!\!\!-\!\!\!-\!\!\!-\!\!\!-\!\!\!\longrightarrow}&
H^{r_1-*}(Y_1){\otimes}{\cdots}{\otimes}H^{r_m-*}(Y_m)\\
\|\wr&&\|\wr\\
H_*(\w X,\w X\setminus Y_1{\times}{\cdots}{\times}Y_m)&\stackrel{\gamma}{\longrightarrow}&
H^{r_1-*}(Y_1{\times}{\cdots}{\times}Y_m)
\end{array}\]
we have $\gamma_{\w\sigma,\omega}\big(a{\otimes}(x_1,\cdots,x_m)\big)=a{\otimes}(\gamma_1(x_1),\cdots,\gamma_m(x_m))$ for all $a\in H^*(K^*_{\w\sigma,\omega})$ and $(x_1,\cdots,x_m)\in T_*^{\w\sigma,\omega}(\w X,M)$, where $\gamma_k$ is as defined in the theorem. So $\zeta_{\sigma,\omega}=\gamma^{-1}_{\w\sigma,\omega}(\partial'_{\sigma,\omega})^{-1}$ is just the isomorphism of the theorem.
\end{proof}\vspace{3mm}

{\bf Example 1.7} For $k\leqslant r$, the standard inclusion of the sphere $S^k\hookrightarrow\RR^{r+1}$ induces an inclusion map $\theta_k\colon S^k\to S^{r+1}$, where we regard $S^{r+1}$ as one-point compactification of $\RR^{r+1}$. $S^{r-k}$ is the deformation retract of the complement space $S^{r+1}{\setminus}\theta_k(S^k)$. For $k_i\leqslant r_i$, $i=1,\cdots,m$, let $M=\ZZ_K\Big(
\begin{array}{ccc}
\scriptstyle{r_1{+}1} &\scriptstyle{\cdots}&\scriptstyle{r_m{+}1}\\
\scriptstyle{ k_1}   &\scriptstyle{\cdots}& \scriptstyle{k_m}\end{array}\Big)=\ZZ_K(X,A)$ be the generalized moment-angle complex with every space pair $(X_i,A_i)=(S^{r_i+1},S^{k_i})$ and the inclusion $\theta_{k_i}\colon S^{k_i}\to S^{r_i+1}$ as defined above. Then $M^c$ is homotopic equivalent to $\ZZ_{K^*}\Big(
\begin{array}{ccc}
\scriptstyle{r_1+\,1\,} &\scriptstyle{\cdots}&\scriptstyle{r_m+\,1\,}\\
\scriptstyle{r_1{-}k_1}   &\scriptstyle{\cdots}& \scriptstyle{r_m{-}k_m}\end{array}\Big)$ and we have $I_M=\{(\sigma,\omega)\,|\,\sigma\in K,\sigma{\cap}\omega=\phi\}$ and $I_{M^c}=\{(\sigma',\omega)\,|\,\sigma'\in K^*,\sigma'{\cap}\omega=\phi\}$. Since all the graded groups $\frak k_{k},\frak c_{k},\frak i_k$ are one dimensional, we use the degree to represent the unique generator of the group. Then by Theorem 1.3, for all $(\sigma,\omega)\in I_M$,
\[H_{k+1}^{\sigma,\omega}\Big(\ZZ_K\Big(
\begin{array}{ccc}
\scriptstyle{r_1{+}1} &\scriptstyle{\cdots}&\scriptstyle{r_m{+}1}\\
\scriptstyle{ k_1}   &\scriptstyle{\cdots}& \scriptstyle{k_m}\end{array}\Big)\Big)=
\oplus_{s+n_1+\cdots n_m=k}\,\w H_s(K_{\sigma,\omega}){\otimes}
(n_1,\cdots,n_m),\]
where $n_i=r_i{+}1$ if $i\in\sigma$, $n_i=0$ if $i\in\w\sigma$, $n_i=k_i$ if $i\in\omega$ and  for all $(\sigma',\omega)\in I_{M^c}$
\[H_{k+1}^{\sigma'\!\!,\omega}\Big(\ZZ_{K^*}\Big(
\begin{array}{ccc}
\scriptstyle{r_1{+}1} &\scriptstyle{\cdots}&\scriptstyle{r_m{+}1}\\
\scriptstyle{r_1{-}k_1}   &\scriptstyle{\cdots}& \scriptstyle{r_m{-}k_m}\end{array}\Big)\Big)
=\oplus_{s+n'_1+\cdots+n'_m=k}\,\w H_s(K^*_{\sigma',\omega}){\otimes}
(n'_1,\cdots,n'_m),\]
where $n'_i=r_i{+}1$ if $i\in\sigma'$, $n'_i=0$ if $i\in\w\sigma'$, $n'_i=r_i{-}k_i$ if $i\in\omega$. Take $\sigma'=\w\sigma$. Notice that $H_s(K_{\sigma,\omega})\cong H^{t-s-1}(K^*_{\w\sigma,\omega})$, where $t$ is the cardinality of $\omega$. The homology and cohomology groups satisfy Theorem 1.6.

\section{ The homology and cohomology of homotopy open covers}\vspace{3mm}

We have computed the (co)homology group of a generalized moment-angle complex, but the coalgebra structure of the homology group (equivalently, the algebra structure of the cohomology group) is not computed. To determine the coalgebra structure, we have to use spectral sequence.\vspace{3mm}

{\bf Definition 2.1} {\it Let $M$ be a topological space. A homotopy open cover $\frak C$ of $M$ is a sequence of topological subspaces $\frak C=(M_1,\cdots,M_n)$ such that every $M_k$ is a deformation retract of an open subspace $M'_k$ of $M$ and that for all $1\leqslant i_1<\cdots<i_s\leqslant n$, $M_{i_1}{\cap}{\cdots}{\cap}M_{i_s}$ is a deformation retract of $M'_{i_1}{\cap}{\cdots}{\cap}M'_{i_s}$ and that $M=\cup_{i=1}^n M'_i$. The cover $\frak C$ is always regarded as an ordered set of symbols, but not a set of subspaces. So $M_i=M_j$ as subspaces for $i\neq j$ is allowed in the sequence. $M_i=\phi$ is also allowed.
}\vspace{3mm}

{\bf Definition 2.2} {\it Let $\frak C=(M_1,\cdots,M_n)$ be a homotopy open cover of $M$. The chain complex $(\frak C_{*,*}(M),d)$ is defined as follows. For $\mu=\{i_0,\cdots,i_s\}\in 2^{[n]}{\setminus}\{\phi\}$, define $M_{\mu}=M_{i_0}{\cap}\cdots{\cap}M_{i_s}$ and $|\mu|=s$. Then
\[\frak C_{*,*}(M)=\oplus_{s,t\geqslant 0}\frak C_{s,t}(M),\quad \frak C_{s,t}(M)=\oplus_{\mu\subset[n],|\mu|=s}\,\{\mu\}{\otimes}H_t(M_{\mu}),\]
where $\{\mu\}$ denotes the $1$-dimensional space generated by $\mu$. For $\mu\subset\nu$, $i_{\nu,\mu}\colon H_*(M_{\nu})\to H_*(M_{\mu})$ is the homomorphism induced by the inclusion map from $M_{\nu}$ to $M_{\mu}$. Then $d\colon \frak C_{s,t}(M)\to \frak C_{s-1,t}(M)$ is defined by
\[d(\{i_0,\cdots,i_s\}{\otimes}a)=\Sigma_{k=0}^s(-1)^{k}\{i_0,\cdots,\wh i_k,\cdots,i_s\}{\otimes}i_{\mu,\mu_k}(a),\quad d(\{i_0\}{\otimes}a)=0,\]
where $\mu=\{i_0,\cdots,i_s\}$, $i_0<\cdots<i_s$, $\mu_k=\{i_0,\cdots,\wh i_k,\cdots,i_s\}$ and $a\in H_*(M_\mu)$.

$(\frak C_{*,*}(M),d)$ has a coproduct $\Delta\colon(\frak C_{*,*}(M),d)\to(\frak C_{*,*}(M){\otimes}\frak C_{*,*}(M),d{\otimes}d)$
defined as follows. For ordered subset $\mu=\{i_0,\cdots,i_s\}$, $\mu'_k=\{i_0,\cdots,i_{k}\}$, $\mu''_k=\{i_k,\cdots,i_{s}\}$ and $a\in H_{|a|}(M_{\mu})$ such that $\Delta_\mu(a)=\Sigma_j a'_j{\otimes}a''_j$, where $\Delta_\mu$ is the coproduct of the homology coalgebra $H_*(M_{\mu})$,
\[\Delta(\{i_0,\cdots,i_s\}{\otimes}a)=\Sigma_j\Sigma_{k=0}^s\,(-1)^{(s-k)|a'_j|}
\{i_0,\cdots,i_k\}{\otimes}i_{\mu,\mu'_k}(a'_j)
\otimes \{i_k,\cdots,i_s\}{\otimes}i_{\mu,\mu''_k}(a''_j).\]
It is easy to check that $(d{\otimes}d)\Delta=\Delta d$  and so the homology $H_{*,*}(\frak C_{*,*}(M),d)$ is a coalgebra over $\kk$.

Dually, the cochain complex $(\frak C^{*,*}(M),\delta)$ is defined as follows.
\[\frak C^{*,*}(M)=\oplus_{s,t\geqslant 0}\frak C^{s,t}(M),\quad \frak C^{s,t}(M)=\oplus_{\mu\subset[n],|\mu|=s}\{\mu\}{\otimes}H^*(M_{\mu}).\]
For $\mu\subset\nu$, $i^*_{\nu,\mu}\colon H^*(M_{\mu})\to H^*(M_{\nu})$ is the homomorphism induced by the inclusion map from $M_{\nu}$ to $M_{\mu}$. Then $\delta\colon \frak C^{s,t}(M)\to \frak C^{s+1,t}(M)$ is defined by that for ordered set $\mu=\{i_0,\cdots,i_k\}$,
\[\delta(\{i_0,\cdots,i_s\}{\otimes}a)=\Sigma_k(-1)^{u+1}
\{i_0,\cdots,i_u,k,i_{u+1}\cdots,i_s\}{\otimes}i^*_{\mu_k,\mu}(a),\quad \delta([n]{\otimes}a)=0,\]
where the sum is taken over all ordered set $\mu_k=\{i_0,\cdots,i_u,k,i_{u+1},\cdots,i_s\}$.

$(\frak C^{*,*}(M),\delta)$ has a  product $\frak C^{*,*}(M){\otimes}\frak C^{*,*}(M)\to \frak C^{*,*}(M)$ defined as follows. $(\mu{\otimes}a)(\nu{\otimes}b)=0$ except that $\mu=\{i_0,\cdots,i_{k}\}$, $\nu=\{i_k,\cdots,i_{s}\}$, $\lambda=\{i_0,\cdots,i_{s}\}$, $i_0<\cdots<i_s$,
\[(\mu{\otimes}a)(\nu{\otimes}b)=(-1)^{|a||\nu|}
\lambda{\otimes}\langle i^*_{\lambda,\mu}(a),i^*_{\lambda,\nu}(b)\rangle_\lambda,\]
where $\langle\,,\rangle_\lambda$ is the product of the cohomology algebra $H^*(M_\lambda)$.
The cohomology $H^{*,*}(\frak C^{*,*}(M),\delta)$ is an algebra over $\kk$.
}\vspace{4mm}

{\bf Theorem 2.3} {\it Let everything be as in Definition 2.2. There is a spectral sequence $(\frak C_{s,t}^r(M),d_r)$ converging to $H_*(M)$ such that every $(\frak C_{*,*}^r(M),d_r)$ is a differential graded coalgebra and there is a differential graded coalgebra isomorphism
\[(\frak C_{*,*}^1(M),d_1)\cong (\frak C_{*,*}(M),d).\]

Dually, there is a spectral sequence $(\frak C^{s,t}_r(M),\delta^r)$ converging to $H^*(M)$ such that every $(\frak C^{*,*}_r(M),\delta^r)$ is a differential graded algebra and there is a differential graded algebra isomorphism\vspace{1mm}\\
\hspace*{61.6mm}$(\frak C^{*,*}_1(M),\delta^1)\cong (\frak C^{*,*}(M),\delta)$.
}

\begin{proof} We only prove the homology case, the cohomology case is just the dual proof. We may suppose that $M_1,\cdots,M_n$ are all open subspaces of $M$.

Define double complex $(U_{*,*},D)$ as follows.
\[U_{*,*}=\oplus_{s,t\geqslant 0}U_{s,t},\quad U_{s,t}=\oplus_{\mu\subset[n],|\mu|=s}\{\mu\}{\otimes}S_t(M_{\mu}),\]
where $S_t(M_{\mu})$ is the singular chain group over $\kk$ freely generated by all singular $t$-simplices of $M_{\mu}$, $\{\mu\}$ and $M_\mu$ are as defined in Definition 2.2. $D\colon T_{s,t}\to T_{s-1,t}{\oplus}T_{s,t-1}$ is defined by that for ordered set $\mu=\{i_0,\cdots,i_s\}\subset[n]$ and singular $t$-simplex $\varpi\in S_t(M_\mu)$,
\[D(\{i_0,\cdots,i_s\}{\otimes}\varpi)=\Sigma_{k=0}^s(-1)^{k}\{i_0,\cdots,\wh i_k,\cdots,i_s\}{\otimes}\varpi+(-1)^s\{i_0,\cdots,i_s\}{\otimes}d(\varpi),\]
where $d$ is the differential of the singular chain complex $(S_*(M_\mu),d)$.

From the horizontal filtration $F_{n}=\oplus_{t\leqslant n}U_{*,t}$ of $(U_{*,*},D)$, we get a spectral sequence $( F_{s,t}^r,d_r)$ converging to $H_*(U_{*,*},D)$. For a singular $t$-simplex $\varpi\in S_t(M_1){+}\cdots{+}S_t(M_m)$, there is a unique subset $\lambda\subset[n]$ such that $\varpi$ is a singular $t$-simplex of $S_t(M_i)$ if and only if $i\in \lambda$. Let $C(\varpi)$ be the subgroup of $U_{*,*}$ generated by all $\mu{\otimes}\varpi$ with $\mu\subset \lambda$. Then $F_{*,*}^1=\oplus_{\varpi}\,H_*(C_*(\varpi),d_1)$, where $\varpi$ is taken over all singular simplices of $S_*(M_1){+}\cdots{+}S_*(M_m)$ and the differential is defined by
\[d_1(\{i_0,\cdots,i_s\}{\otimes}\varpi)=\Sigma_{k=0}^s(-1)^{k}\{i_0,\cdots,\wh i_k,\cdots,i_s\}{\otimes}\varpi.\]
It is obvious that $(C_*(\varpi),d_1)\cong(C_*(2^\lambda),d)$, where $(C_*(2^{\lambda}),d)$ is the simplicial chain complex over $\kk$ of the full simplicial complex $2^\lambda$ and so $H_0(C_*(\varpi),d_1)=H_0(2^{\lambda})=\kk$ and $H_k(C_*(\varpi),d_1)=H_k(2^{\lambda})=0$ if $k\neq 0$. The generator of $H_0(C_*(\varpi),d_1)$ is represented by $\{i\}{\otimes}\varpi$ for any $i\in\lambda$, i.e., $\{i\}{\otimes}\varpi$ and $\{j\}{\otimes}\varpi$ represent the same homology class in $F^1_{*,*}$ for all $i,j\in\lambda$. This implies that the correspondence $\{i\}{\otimes}\varpi\to\varpi$ induces a chain complex isomorphism \[(F^1_{*,*},d_1)=(F^1_{0,*},d_1)\cong(S_*(M_1){+}\cdots{+}S_*(M_m),d).\]
Since $d_2\colon F^2_{s,t}\to F^2_{s+1,t-2}$, the spectral sequence collapse from $r\geqslant 2$ and  $H_{0,*}(F^1_{*,*},d_1)\cong H_{*}(U_{*,*},D)$. By excision axiom, $(S_*(M_1){+}\cdots{+}S_*(M_m),d)$ is homotopic equivalent to the singular chain complex $(S_*(M),d)$. So $H_{0,*}(F^1_{*,*},d_1)\cong H_*(M)$.

From the vertical filtration $E_{n}=\oplus_{s\leqslant n}U_{s,*}$ of $(U_{*,*},D)$, we get a spectral sequence $(E_{s,t}^r,d_r)$converging to $H_*(U_{*,*},D)=H_*(M)$. Define $(E_{s,t}^r,d_r)=(\frak C_{s,t}^r(M),d_r)$.  By definition, $(\frak C_{*,*}^1(M),d_1)=(\frak C_{*,*}(M),d)$.

Define coproduct $\w\Delta\colon (U_{*,*},D)\to (U_{*,*},D){\otimes}(U_{*,*},D)$ as follows. For a singular simplex $\varpi\in S_t(M_{\mu})$ with $\mu=\{i_0,\cdots,i_s\}$ and $\Delta_\mu(\varpi)=\varpi'{\otimes}\varpi''$ (abbreviation of $\Sigma_j\varpi'_j{\otimes}\varpi''_j$), where $\Delta_\mu\colon S_*(M_\mu)\to S_*(M_\mu){\otimes}S_*(M_\mu)$ is the homomorphism of singular chain groups induced by the diagonal map of $M_{\mu}$, define
\[\w\Delta(\{i_0,\cdots,i_s\}{\otimes}\varpi)=\Sigma_{k=0}^s(-1)^{(s-k)|\varpi'|}
\{i_0,\cdots,i_k\}{\otimes}\varpi'\otimes\{i_k,\cdots,i_s\}{\otimes}\varpi''.\]
Then\\
\hspace*{20mm}$\quad(D{\otimes}D)(\w \Delta(\{i_0,\cdots,i_s\}{\otimes}\varpi))$\\
\hspace*{20mm}$=(D{\otimes}D)((-1)^{(s-k)|\varpi'|}
\{i_0,\cdots,i_k\}{\otimes}\varpi'\otimes\{i_k,\cdots,i_s\}{\otimes}\varpi'')$\\
\hspace*{20mm}$=\,\,\,\,\Sigma_{k=0}^s\Sigma_{u=0}^k(-1)^{(s-k)|\varpi'|+u}
\{i_0,\cdots,\wh i_{u},\cdots,i_k\}{\otimes}\varpi'\otimes\{i_k,\cdots,i_s\}{\otimes}\varpi''$\\
\hspace*{20mm}$\quad{+}\Sigma_{k=0}^s(-1)^{(s-k)|\varpi'|+k}
\{i_0,\cdots,i_k\}{\otimes}d\varpi'\otimes\{i_k,\cdots,i_s\}{\otimes}\varpi''$\\
\hspace*{20mm}$\quad{+}\Sigma_{k=0}^s\Sigma_{u=0}^{s-k}(-1)^{(s-k)|\varpi'|+k+|\varpi'|+u}
\{i_0,\cdots,i_k\}{\otimes}\varpi'\otimes\{i_k,\cdots,\wh i_{k+u},\cdots,i_s\}{\otimes}\varpi''$\\
\hspace*{20mm}$\quad+\Sigma_{k=0}^s(-1)^{(s-k)|\varpi'|+s+|\varpi'|}
\{i_0,\cdots,i_k\}{\otimes}\varpi'\otimes\{i_k,\cdots,i_s\}{\otimes}d\varpi''$\\
\hspace*{20mm}$=\,\,\,\,\Sigma_{k=0}^s\Sigma_{u=0}^{k-1}(-1)^{(s-k)|\varpi'|+u}
\{i_0,\cdots,\wh i_{u},\cdots,i_k\}{\otimes}\varpi'\otimes\{i_k,\cdots,i_s\}{\otimes}\varpi''$\\
\hspace*{20mm}$\quad{+}\Sigma_{k=0}^s\Sigma_{u=1}^{s-k}(-1)^{(s-k)|\varpi'|+k+|\varpi'|+u}
\{i_0,\cdots,i_k\}{\otimes}\varpi'\otimes\{i_k,\cdots,\wh i_{k+u},\cdots,i_s\}{\otimes}\varpi''$\\
\hspace*{20mm}$\quad{+}\Sigma_{k=0}^{s-1}(-1)^{(s-k)|\varpi'|+k}
\{i_0,\cdots,\wh i_{k}\}{\otimes}\varpi'\otimes\{i_k,\cdots,i_s\}{\otimes}\varpi''$\\
\hspace*{20mm}$\quad{+}\Sigma_{k=1}^{s}(-1)^{(s-k)|\varpi'|+k+|\varpi'|}
\{i_0,\cdots,i_k\}{\otimes}\varpi'\otimes\varpi''{\otimes}\{\wh i_{k},\cdots,i_s\}$\\
\hspace*{20mm}$\quad{+}\Sigma_{k=0}^s(-1)^{(s-k)(|\varpi'|-1)+s}
\{i_0,\cdots,i_k\}{\otimes}d\varpi'\otimes\{i_k,\cdots,i_s\}{\otimes}\varpi''$\\
\hspace*{20mm}$\quad{+}\Sigma_{k=0}^s(-1)^{(s-k)|\varpi'|+s+|\varpi'|}
\{i_0,\cdots,i_k\}{\otimes}\varpi'\otimes\{i_k,\cdots,i_s\}{\otimes}d\varpi''$\\
\hspace*{20mm}$=\w\Delta(\Sigma_{u=0}^s(-1)^{u}\{i_0,\cdots,\wh i_u,\cdots,i_s\}{\otimes}\varpi)
{+}\w\Delta((-1)^s\{i_0,\cdots,i_s\}{\otimes}d\varpi)$\\
\hspace*{20mm}$=\w\Delta(D(\{i_0,\cdots,i_s\}{\otimes}\varpi))$\\

Thus, $\w\Delta$ induces spectral sequence homomorphism $\w\Delta_r\colon
(E^r_{*,*},d_r)\to(E^r_{*,*}{\otimes}E^r_{*,*},d_r{\otimes}d_r)$. By K${\rm \ddot{u}}$nneth Theorem, $H_*(U_{*,*}{\otimes}U_{*,*},D{\otimes}D)\cong H_*(M){\otimes}H_*(M)\cong H_*(M{\times}M)$.

$\w\Delta$ also induces a spectral sequence homomorphism $\o\Delta_r\colon
(F^r_{*,*},d_r)\to(F^r_{*,*}{\otimes}F^r_{*,*},d_r{\otimes}d_r)$. By definition, for a homology class in $F^1_{0,*}$ represented in $U_{*,*}$ by $\{i\}{\otimes}[\varpi]$ such that $\Delta_i(\lambda)=\lambda'{\otimes}\lambda''$ ($\Delta_i$ is the coproduct of $H_*(M_i)$),  $\o\Delta_1(\{i\}{\otimes}\lambda)=\{i\}{\otimes}\lambda'\otimes\{i\}{\otimes}\lambda''$. This implies the following diagram commutes.
$$
\begin{array}{ccc}
(F^1_{0,*},d_1) &\stackrel{\o\Delta_1}{\longrightarrow}&(F^1_{0,*}{\otimes}F^1_{0,*},d_1{\otimes}d_1) \\
\|\wr &   & \|\wr \\
(S_*(M),d) &\stackrel{\Delta}{\longrightarrow}&(S_*(M){\otimes}S_*(M),d)
\end{array}
$$
where $\Delta$ is the homomorphism of singular chain groups induced by the diagonal map of $M$. Since the spectral sequence $F^{r}_{s,t}$ and $F^{r}_{s,t}{\otimes}F^{r}_{s,t}$ collapse from $r\geqslant 2$,  we have the following commutative diagram
\[
\begin{array}{ccc}
H_*(U_{*,*},D) &\stackrel{\w\Delta_*}{\longrightarrow}&H_*(U_{*,*}{\otimes}U_{*,*},D{\otimes}D) \\
\|\wr &   & \|\wr \\
H_*(M) &\stackrel{\Delta_*}{\longrightarrow}&H_*(M{\times}M).
\end{array}\]
So $H_*(U_{*,*},D)\cong H_*(S_*(M),d)$ is a coalgebra isomorphism.
\end{proof}

Notice that if we do not consider the (co)algebra strucure, the above theorem holds for (co)homology over an Abellian group $G$. But all other main conclusions (such as direct sum decomposition, generating set) in this paper can not be generalized to the case of an Abellian group $G$. So we only discuss field case.\vspace{3mm}

{\bf Definition 2.4} {\it Let everything be as in Definition 2.2. For $x=\mu{\otimes}a\in \frak C_{*,*}(M)$, the support complex $(\frak C_{*,*}(x),d)$ of $x$ is the chain subcomplex of $(\frak C_{*,*}(M),d)$ generated by all $\lambda{\otimes}i_{\mu,\lambda}(a)$ such that $\lambda\subset\mu$. Define $\frak C_x=\{\lambda\subset\mu\,|\,i_{\mu,\lambda}(a)=0\}$. It is obvious that $\frak C_x$ is a simplicial complex.
}\vspace{3mm}

{\bf Theorem 2.5} {\it Let everything be as in Definition 2.4. For $x=\mu{\otimes}a\in \frak C_{|\mu|,|a|}(M)$, there is a chain complex isomorphism $(\frak C_{*,|a|}(x),d)=(C_*(2^\mu,\frak C_x),d)$, where $C_*(2^\mu,\frak C_x),d)$ is the relative simplicial chain complex over $\kk$, i.e., the quotient complex of $(C_*(2^\mu),d)$ over $(C_*(\frak C_x),d)$. So we may define $H_s(\frak C_{*,*}(x),d)=H_{s,|a|}(\frak C_{*,|a|}(x),d)\cong H_s(2^\mu,\frak C_x)\cong \w H_{s-1}(\frak C_x)$.
}
\begin{proof} The correspondence $\lambda{\otimes}i_{\mu,\lambda}(a)\to\lambda$ for all $i_{\mu,\lambda}(a)\neq 0$ induces a chain complex isomorphism from $(\frak C_{*,*}(x),d)$ to $(C_*(2^\mu,\frak C_x),d)$.
\end{proof}

{\bf Definition 2.6} {\it Let $\frak C$ be a homotopy open cover of $M$. $\frak C$ is a simple cover
(with respect to $\kk$) if there is a subset $T_{*,*}(\frak C)$ of $\frak C_{*,*}(M)$ such that
\[(\frak C_{*,*}(M),d)=\oplus\,_{x\in T_{*,*}(\frak C)}(\frak C_{*,*}(x),d)\]
and the spectral sequence $\frak C^r_{s,t}(M)$ defined in Theorem 2.3 collapse from $r\geqslant 2$ and so for all $k$,
\[H_k(X)=\oplus_{s+t=k}\frak C_{s,t}^2(M)=\oplus_{s+t=k}\left(\oplus\,_{x\in T_{*,t}(\frak C)}\,\,H_{s}(\frak C_{*,*}(x),d)\right).\]
$T_{*,*}(\frak C)$ is called a generating set of $\frak C_{*,*}(M)$.
}\vspace{3mm}
\section {Homology coalgebra and cohomology algebra}\vspace{3mm}

{\bf Definition 3.1} {\it Let $K$ be a simplicial complex. A simplex cover of $K$ is a sequence $\frak S=(\sigma_1,\cdots,\sigma_n)$ of simplices (repetition and empty simplex allowed) of $K$ that contains the set of maximal simplices of $K$, i.e., every simplex of $K$ is the face of some $\sigma_k$.

For a generalized moment-angle complex $M=\ZZ_K(X,A)$, $(X,A)=\{(X_k,A_k)\}_{k=1}^m$, and  a simplex cover $\frak S=(\sigma_1,\cdots,\sigma_n)$ of $K$, the sequence $\frak C=(M_1,\cdots,M_n)$ with $M_k=D(\sigma_k)$ ($D(\sigma_k)$ as defined in Definition 1.1) is the cover of $M$ associated to $\frak S$.}\vspace{3mm}

{\bf Definition 3.2} {\it Let $M=\ZZ_K(X,A)$, $(X,A)=\{(X_k,A_k)\}_{k=1}^m$, be a generalized moment-angle complex. Then for a simplex cover $\frak S=(\sigma_1,\cdots,\sigma_n)$ of $K$, the associated cover $\frak C=(M_1,\cdots,M_n)$ is a homotopy open cover of $M$ and by Definition 2.2, there is a chain complex $(\frak C_{*,*}(M),d)$.

The base $\Gamma_{*,*}(\frak C)$ of $\frak C_{*,*}(M)$ is taken as follows. Take a fixed base for $\frak k_k,\frak c_k,\frak i_k$ (as defined in Definition 1.2). By K\"{u}nneth Theorem, for a space $Y=Y_1{\times}{\cdots}{\times}Y_m$ with $Y_k=A_k$ or $X_k$ for every $k$, the base of $H_*(Y)$ is the set of all $\xx=(x_1,{\cdots},x_m)$ with every $x_k$ a base elements of $H_*(Y_k)$, where we do not distinguish a base element in ${\rm coim}i_k$ and its image in ${\rm im}i_k$. Then $\Gamma_{*,*}(\frak C)$ is the set of all $\mu{\otimes}(x_1,\cdots,x_m)$ such that $\mu=\{i_0,\cdots,i_s\}\subset[n]$ and that $(x_1,\cdots,x_m)$ is a base element of $H_*(M_\mu)$, where $M_\mu=D(\sigma_{i_0}){\cap}{\cdots}{\cap}D(\sigma_{i_s})=D({\cap}_{j\in\mu}\sigma_j)=Y_1{\times}{\cdots}{\times}Y_m$ with $Y_k=X_k$ if $k\in\cap_{j\in\mu}\sigma_j$ and $Y_k=A_k$ otherwise. Equivalently, $\Gamma_{*,*}(\frak C)$ is the following set
\[\{\mu{\otimes}(x_1,\cdots,x_m)\,|\,\mu\subset[n],\,(x_1,\cdots,x_m)\in T_*(M),\,
x_k\in\frak c_k{\cup}\frak i_k,\,{\rm if}\,k\in\cap_{j\in\mu}\sigma_j,\,
x_k\in\frak k_k{\cup}\frak i_k,\,{\rm if}\,k\not\in\cap_{j\in\mu}\sigma_j\},\]
where we use the same symbol to denote a vector space and its base and so $T_*(M)=\cup_{(\sigma,\omega)\in I_M}T_*^{\sigma,\omega}(M)$ with $T_*^{\sigma,\omega}(M)$ as defined in Definition 1.2.
}\vspace{3mm}

{\bf Theorem 3.3} {\it Let everything be as in Definition 3.2. Give $\Gamma_{*,*}(\frak C)$ a graded graph structure defined as follows. For $\mu{\otimes}\xx\in \Gamma_{*,*}(\frak C)$, its degree is $|\mu|$, the cardinality of $\mu$ minus $1$. For $\mu{\otimes}\xx,\nu{\otimes}\yy\in \Gamma_{*,*}(\frak C)$ with $\mu$ a proper subset of $\nu$, they are the two vertices of an edge if and only if $i_{\nu,\mu}(\yy)=\xx$ and $|\mu|=|\nu|{-}1$, where $i_{\nu,\mu}$ is as defined in Definition 2.2. Then the graded graph $\Gamma_{*,*}(\frak C)$ is simple in the sense that every connected component of it has but one vertex with maximal degree which is called the top vertex of the component. So the connected component decomposition of the graph $\Gamma_{*,*}(\frak C)$ corresponds to the direct sum decomposition of the chain complex $(\frak C_{*,*}(M),d)$. Precisely, let $T_{*,*}(\frak C)$ be the set of all top vertices of $\Gamma_{*,*}(\frak C)$, then $\frak C$ is a simple cover (as in Definition 2.6) of $\ZZ_K(X,A)$ with generating set  $T_{*,*}(\frak C)$ and direct sum decomposition
\[H_k(M)=\oplus_{s+t=k}H_{s,t}(\frak C_{*,*}(M),d)=\oplus_{s+t=k}\left(\oplus\,_{x\in T_{*,t}(\frak C)}\,\,H_{s}(\frak C_{*,*}(x),d)\right),\]
where $(\frak C_{*,*}(x),d)$ is as defined in Definition 2.4 and $H_*(\frak C_{*,*}(x),d)$ is as defined in Theorem 2.5.}
\begin{proof} Notice that since we do not distinguish coim$i_k$ and im$i_k$, the homomorphism $i_k$ restricted on coim is the identity isomorphism. So for $\lambda{\otimes}\xx,\mu{\otimes}\yy\in\Gamma_{*,*}(\frak C)$, if $i_{\mu,\lambda}(\yy)=\xx$, then $\xx=\yy$. This implies that all connected components of $\Gamma_{*,*}(\frak C)$ is a graph with vertex set of the form $\{\mu_1{\otimes}\xx,\cdots,\mu_s{\otimes}\xx\}$ with $\xx\in T_*(M)$, i.e., the set of connected components of $\Gamma_{*,*}(\frak C)$ is in 1-1 correspondence with $T_*(M)$. For a connected component with vertex set $\{\mu_1{\otimes}\xx,\cdots,\mu_s{\otimes}\xx\}$, $x=(\cup_{i=1}^s\mu_i){\otimes}\xx$ is the unique top vertex of the component and the connected component chain complex is obviously $(C_{*,*}(x),d)$.

Now we prove the spectral sequence $E^r_{s,t}=\frak C^r_{s,t}(M)$ as defined in Theorem 2.3 collapses from $r\geqslant 2$. Let $(U_{*,*},D)$ and everything else be as in the proof of Theorem 2.3.  Then $U_{*,*}$ is generated by elements of the form $\mu{\otimes}(x_1,\cdots,x_m)$ such that $x_k\in S_*(X_k)$ if $k\in\cap_{j\in\mu}\sigma_j$ and $x_k\in S_*(A_k)$ if $k\not\in\cap_{j\in\mu}\sigma_j$. Let $(C_*(M),d)$ be as defined in the proof of Theorem 1.3. We use the same symbol to denote both the homology class and a representative of it, both the space and a base of it. Then $(C_*(M),d)$ is a chain subcomplex of $(S_*(X_1){\otimes}{\cdots}{\otimes}S_*(X_m),d{\otimes}{\cdots}{\otimes}d)$ by taking representatives as in the proof of Theorem 1.3. So we may define double subcomplex $(\w {\frak C}_{*,*},D)$ of $(U_{*,*},D)$ generated by all $\mu{\otimes}(x_1,\cdots,x_m)$ such that $x_k\in \frak k_k{\cup}\frak i_k{\cup}\frak c_k{\cup}\frak q_k$ if $k\in\cap_{j\in\mu}\sigma_j$ and $x_k\in \frak k_k{\cup}\frak i_k$ if $k\not\in\cap_{j\in\mu}\sigma_j$. Then by Definition 3.2, $\Gamma_{*,*}(\frak C)$ is a subgroup (not a subcomplex) of $\w{\frak C}_{*,*}$ and $(\w{\frak C}_{*,*},D)$ is the smallest subcomplex of $(U_{*,*},D)$ containing $\Gamma_{*,*}(\frak C)$.

Let $(Z_{*,*},D)$ be the subcomplex of $(\w{\frak C}_{*,*},D)$ generated by all $\mu{\otimes}(x_1,{\cdots},x_m)$ such that there is at least one $k\in\cap_{j\in\mu}\sigma_j$ such that $x_k\in \frak k_k$ or $\frak q_k$.
Let $(Z'_{*,*},D)$ be the subcomplex of $(Z_{*,*},D)$ generated by all $\mu{\otimes}(x_1,{\cdots},x_m)$ such that there is no $x_k\in\frak q_k$. Let $(F_n,D)$ be the subcomplex of $(Z_{*,*}(x),D)$ generated by all $\lambda{\otimes}\yy$ such that $|\lambda|\leqslant n$. Then $\{F_n\}$ is a filtration that induces a spectral sequence $(F_{s,t}^r,d_r)$ converging to $H_*(Z_{*,*},D)$. For $\mu{\otimes}(x_1,{\cdots},x_m)\in Z'_{*,*}$, define $(D(x_k),d'_k)$ to be the acylic chain complex such that $D(x_k)$ is generated by $\{x_k,x'_k\}$ with $d'_k(x'_k)=x_k$ if $x_k\in\frak k_k$ and $k\in\cap_{j\in\mu}\sigma_j$ and define $(D(x_k),d'_k)$ to be the trivial chain complex with $D(x_k)$ generated by $x_k$ if  $x_k\not\in\frak k_k$ or $k\not\in\cap_{j\in\mu}\sigma_j$. By definition, $(F_{*,*}^1,d_1)\cong\oplus_{\mu{\otimes}(x_1,{\cdots},x_m)\in Z'_{*,*}}(D(x_1){\otimes}{\cdots}{\otimes}D(x_m),d'_1{\otimes}{\cdots}{\otimes}d'_m)$ with the isomorphism induced by the correspondence $\mu{\otimes}\xx\to \xx$. So $F_{*,*}^2=0$ and $H_*(Z_{*,*},D)=0$. It is obvoious that the quotient complex of $(\w{\frak C}_{*,*},D)$ over $(Z_{*,*},D)$ is isomorphic to $(\frak C_{*,*}(M),d)$. So $H_*(\w{\frak C}_{*,*},D)=H_*(\frak C_{*,*}(M),d)$. Since $(\w{\frak C}_{*,*},D)$ is a double subcomplex of $(U_{*,*},D)$, all the homology classes of $H_*(\frak C_{*,*}(M),d)=\frak C^2_{*,*}(M)$ survive to infinity in the spectral sequence.
\end{proof}

{\bf Definition 3.4} {\it Let $M=\ZZ_K(X,A)$, $(X,A)=\{(X_k,A_k)\}_{k=1}^m$ be a generalized moment-angle complex and $\frak C$ be a cover associated to the simplex cover $\frak S=(\sigma_1,\cdots,\sigma_n)$ of $K$. For $(\sigma,\omega)\in I_M$ ($I_M$ as in Definition 1.2), the simplicial complex $\frak S_{\sigma,\omega}$ is defined as follows. The vertex set of $\frak S_{\sigma,\omega}$ is the subset $[\sigma,\omega]=\{k\in[n]\,|\,\sigma\subset\sigma_k,\,\sigma_k{\cap}\omega\neq\phi\}$ and $\{i_0,\cdots,i_s\}\subset [\sigma,\omega]$ is a simplex of $\frak S_{\sigma,\omega}$ if and only if $(\cap_{j=0}^s\sigma_{i_j}){\cap}\omega\in K$.

For $x=\mu{\otimes}(x_1,\cdots,x_m)\in \Gamma_{*,*}(\frak C)$, define $\rho(x)=\{k\,|\,x_k\in\frak c_k\}$ and $\varrho(x)=\{k\,|\,x_k\in\frak k_k\}$ and $T_{*,*}^{\sigma,\omega}(\frak C)=\{x\in T_{*,*}(\frak C)\,|\,\rho(x)=\sigma,\,\varrho(x)=\omega\}$.
}\vspace{3mm}

Notice the difference between $K_{\sigma,\omega}$ in Definition 1.2 and $\frak S_{\sigma,\omega}$ in the above definition. The vertex set $\langle\sigma,\omega\rangle$ of $K_{\sigma,\omega}$ is a subset of $[m]$ and the vertex set $[\sigma,\omega]$ of $\frak S_{\sigma,\omega}$ is a subset of $[n]$. But $[m]$ and $[n]$ are two irrelevant sets, since $[m]$ is a substitute of $(X,A)$ and $[n]$ is a substitute of $\frak S$. We always use $\lambda,\mu,\nu,\mu_i,\mu',\cdots$ to denote subsets of $[n]$ and other Greek letters to denote subsets of $[m]$. \vspace{3mm}

{\bf Theorem 3.5} {\it Let everything be as in Theorem 3.3 and Definition 3.4. Then the following conclusions hold.

1) For any $x=\mu{\otimes}(x_1,\cdots,x_m)\in \Gamma_{*,*}(\frak C)$ with $\rho(x)=\sigma$, $\varrho(x)=\omega$ and $\frak C_x$ as defined in Definition 2.4,  $\frak C_x=\frak S_{\sigma,\omega}|_\mu=\{\lambda{\cap}\mu\,|\,\lambda\in\frak S_{\sigma,\omega}\}$ and there is a chain complex isomorphism $(\frak C_{*,*}(x),d)\cong(C_*(2^{\mu},\frak S_{\sigma,\omega}|_\mu),d)$. Specifically, if $x\in T_{*,*}(\frak C)$, $(C_{*,*}(x),d)\cong(C_*(2^\mu,\frak S_{\sigma,\omega}),d)$.

2) $\mu{\otimes}\xx\to\xx$ is a 1-1 correspondence from $T_{*,*}^{\sigma,\omega}(\frak C)$ to the base of $T_{*}^{\sigma,\omega}(M)$ for all $(\sigma,\omega)\in I_M$, where $T_{*}^{\sigma,\omega}(M)$ is as defined in Definition 1.2.

3) The geometrical realization space $|\frak S_{\sigma,\omega}|\simeq|K_{\sigma,\omega}|$ ($\simeq$ means homotopic equivalent of topological spaces), where $K_{\sigma,\omega}$ is as in Definition 1.2. So for $x=\mu{\otimes}\xx\in T_{*,*}^{\sigma,\omega}(\frak C)$,\\
\hspace*{20mm}$H_k(\frak C_{*,*}(x),d)=H_k(2^{\mu},\frak S_{\sigma,\omega})\cong\w H_{k-1}(\frak S_{\sigma,\omega})\cong\w H_{k-1}(K_{\sigma,\omega})\cong H_k(2^{\langle \sigma,\omega\rangle},K_{\sigma,\omega}).$
}

\begin{proof} 1) The correspondence $\lambda{\otimes}(x_1,\cdots,x_m)\to \lambda$ induces a chain complex isomorphism from $(\frak C_{*,*}(x),d)$ to $(C_*(2^{\mu},\frak S_{\sigma,\omega}|_\mu),d)$.

2) Just the conclusion that $T_*(M)$ is in 1-1 correspondence with the set of top vertices $T_{*,*}(\frak C)$ of $\Gamma_{*,*}(\frak C)$ which is proven in Theorem 3.3.

3) Suppose $\sigma=\{i_0,\cdots,i_s\}$. Then $\frak S|_\omega=(\sigma_{i_0}{\cap}\omega,\cdots,\sigma_{i_s}{\cap}\omega)$ is a simplex cover of $K_{\sigma,\omega}$. Let $U_k=|{\rm star}_{K_{\sigma,\omega}}(\sigma_{i_k}{\cap}\omega)|{\setminus}|{\rm link}_{K_{\sigma,\omega}}(\sigma_{i_k}{\cap}\omega)|$. Then $\{U_0,\cdots,U_s\}$ is an open cover of $|K_{\sigma,\omega}|$ such that every non-empty $U_{j_0}{\cap}\cdots{\cap}U_{j_t}$ is contractible. By definition, $\frak S_{\sigma,\omega}$ is the nerve of the open cover $\{U_0,\cdots,U_s\}$. So $|\frak S_{\sigma,\omega}|\simeq|K_{\sigma,\omega}|$.
\end{proof}

{\bf Theorem 3.6} {\it Let everything be as in Theorem 3.5. For $x_i\in T_{*}^{\sigma_i,\omega_i}(\frak C)$, $i=1,2,3$, such that $\sigma_2\subset\sigma_1$, $\sigma_3\subset\sigma_1$ and $\omega_1\subset\omega_2{\cup}\omega_3$, there is a chain complex homomorphism
\[\w\Delta_{x_1}^{x_2,x_3}\colon (\frak C_{*,*}(x_1),d)\to
(\frak C_{*,*}(x_2){\otimes}\frak C_{*,*}(x_3),d{\otimes}d)\]
that induces a homology homomorphism $\Delta_{x_1}^{x_2,x_3}$ and the following commutative diagram ($v_i=\langle\sigma_i,\omega_i\rangle$)
\[\begin{array}{ccc}
H_*(\frak C_{*,*}(x_1),d)&\stackrel{\Delta_{x_1}^{x_2,x_3}}{-\!\!\!-\!\!\!-\!\!\!-\!\!\!-\!\!\!\longrightarrow}
&H_*(\frak C_{*,*}(x_2),d){\otimes}H_*(\frak C_{*,*}(x_3),d)\\
    \|\wr &&\|\wr\\
H_{*}(|2^{v_1}|,|K_{\sigma_1,\omega_1}|)&\stackrel{\Delta_{v_1}^{v_2,v_3}}{-\!\!\!-\!\!\!-\!\!\!-\!\!\!-\!\!\!\longrightarrow}
&\,\,\,\,H_{*}(|2^{v_2}|,|K_{\sigma_2,\omega_2}|){\otimes}H_{*}(|2^{v_3}|,|K_{\sigma_3,\omega_3}|),
\end{array}\]
where the vertical isomorphisms are obtained by replacing the simplicial homology on the right side of the isomorphism in 3) of Theorem 3.5 by singular homology and $\Delta_{v_1}^{v_2,v_3}$ is the coproduct between relative singular homology groups defined as follows. For simplicial complexes $K,K',K''$ with respectively vertex set $\varrho,\varrho',\varrho''$ such that $K$ is a simplicial subcomplex of $K'{\cup}K''$ and $S$ a set such that $\varrho'{\cup}\varrho''\!\subset\! S$, the coproduct $\Delta_\varrho^{\varrho',\varrho''}$ is defined by the following commutative diagram, which is obviously independent of the choice of $S$
\[\begin{array}{ccc}
H_{*}(|2^{\varrho}|,|K|)&\stackrel{\Delta_\varrho^{\varrho',\varrho''}}{-\!\!\!-\!\!\!-\!\!\!-\!\!\!-\!\!\!\longrightarrow}
&H_{*}(|2^{\varrho'}|,|K'|){\otimes}H_{*}(|2^{\varrho''}|,|K''|)\\
\|\wr&&\|\wr\\
H_{*}(|2^{S}|,|K|)&\stackrel{\Delta_*}{-\!\!\!-\!\!\!-\!\!\!-\!\!\!-\!\!\!\longrightarrow}
&H_{*}(|2^{S}|,|K'|){\otimes}H_{*}(|2^{S}|,|K''|),\\
\end{array}\]
where the vetical isomorphisms are induced by the inclusion map and $\Delta(x)=(x,x)$ for all $x\in|2^S|$ is the diagonal map of $|2^{S}|$ regarded as a map between space pairs from the space pair $(|2^{S}|,|K|)$ to the product space pair $(|2^{S}|,|K'|){\times}(|2^{S}|,|K''|)=
(|2^{S}|{\times}|2^{S}|,|2^{S}|{\times}|K''|\cup|K'|{\times}|2^{S}|)$.

}
\begin{proof} For the simplex cover $\frak S=(\sigma_1,\cdots,\sigma_n)$ of $K$, define simplicial complex $\o K$ as follows. The vertex set of $\o K$ is $[n]$. An non-empty set $\{i_0,{\cdots},i_s\}\subset[n]$ is a simplex of $\o K$ if and only if $\cap_{j=0}^s\sigma_{i_j}\neq\phi$. It is obvious that $\o K$ is a simplicial complex. Define maps $\psi\colon K{\setminus}\{\phi\}\to \o K{\setminus}\{\phi\}$ and $\varphi\colon \o K{\setminus}\{\phi\}\to K{\setminus}\{\phi\}$ as follows. For a non-empty simplex $\sigma$ of $K$, $\psi(\sigma)=\{k\in[n]\,|\,\sigma{\subset}\sigma_{k}\}$. For a non-empty simplex $\mu=\{i_0,{\cdots},i_s\}$ of $\o K$, $\varphi(\mu)={\cap}_{j=0}^s\sigma_{i_j}$. $\psi$ and $\varphi$ are order reversing, i.e., for $\sigma\subset\tau$, $\psi(\sigma)\supset\psi(\tau)$ and for $\mu\subset\nu$, $\varphi(\mu)\supset\varphi(\nu)$. So $\psi$ and $\varphi$ are not simplicial maps but they induce simplicial maps between barycentric subdivision complexes.

For a simplicial complex $L\neq\{\}$, denote the barycentric subdivision complex of $L$ by $S(L)$. $S(L)=L$ if $L=\{\phi\}$. Then $S(L){\setminus}\{\phi\}$ is in 1-1 correspondence with the set of all set of non-empty simplices $\{\sigma_1,\cdots,\sigma_s\}$ of $L$ such that $\sigma_i$ is a proper face of $\sigma_{i+1}$ for $i=1,\cdots,s{-}1$. Define $\Psi\colon S(2^{K{\setminus}\{\phi\}})\to S(2^{\o K{\setminus}\{\phi\}})$ and $\Phi\colon S(2^{\o K{\setminus}\{\phi\}})\to S(2^{K{\setminus}\{\phi\}})$ as follows. For a non-empty simplex $\{\eta_1,{\cdots},\eta_s\}$ of $S(2^{\o K{\setminus}\{\phi\}})$, $\Psi(\{\eta_1,{\cdots},\eta_s\})=\{\psi(\eta_s),{\cdots},\psi(\eta_1)\}$ and for a non-empty simplex $\{\mu_1,{\cdots},\mu_t\}$ of $S(2^{\o K{\setminus}\{\phi\}})$, $\Phi(\{\mu_1,{\cdots},\mu_t\})=\{\varphi(\mu_t),{\cdots},\varphi(\mu_1)\}$, where repetition is canceled in the image set. It is obvious that $\Psi$ and $\Phi$ are simplicial maps.
For a simplicial map $f\colon K\to L$, denote the corresponding piecewise linear map between the geometrical realization spaces by $|f|\colon |K|\to |L|$. Notice that for a non-empty simplex $\sigma\in K$, $\sigma\subset\phi\psi(\sigma)$. If $\sigma\neq\phi\psi(\sigma)$, then for all $\eta\subset\sigma$, $\phi\psi(\eta)=\phi\psi(\sigma)$. This implies that if $\sigma=\phi\psi(\sigma)$, then for all $\sigma\subset\tau$, $\tau=\phi\psi(\tau)$. Thus, for any simplex $V=\{\eta_1,\cdots,\eta_s\}$ of $S(2^{K\setminus\{\phi\}})$, if $\Phi\Psi(V)\neq V$, then there is a $0\leqslant k\leqslant s$ such that $\eta_k\neq\phi\psi(\eta_k)$ and $\eta_{j}=\phi\psi(\eta_{j})$ for $j>k$. So both $V$ and $\Phi\Psi(V)$ are faces of the same simplex $\{\eta_1,\cdots,\eta_k,\phi\psi(\eta_{k}),\eta_{k+1},\cdots,\eta_s\}$ (cancel $\phi\psi(\eta_k)$ if $\phi\psi(\eta_k)=\eta_{k+1}$) of $S(2^{K\setminus\{\phi\}})$. This implies that $|\Phi||\Psi|$ is homotopic to the identity map of $|S(2^{K\setminus\{\phi\}})|$ by the linear homotopy $H(t,x)=(1{-}t)\Phi\Psi(x)+tx$ for all $x\in|S(2^{K\setminus\{\phi\}})|$. For a simplex $W=\{\mu_{1},{\cdots},\mu_{r_1},\cdots,\mu_{r_{t-1}+1},{\cdots},\mu_{r_t}\}$ of $S(2^{\o K\setminus\{\phi\}})$ such that $\cap_{k\in\mu_{s}}\sigma_k=\sigma_{u_i}$ for $s=r_{i-1}{+}1,\cdots,r_i$ and every $\sigma_{u_{i+1}}$ is a proper subset of $\sigma_{u_{i}}$, both $W$ and $\Psi\Phi(W)$ are faces of the same simplex $\{\mu_{1},{\cdots},\mu_{r_1},\psi\phi(\mu_{r_1}),\cdots,\mu_{r_{t-1}+1},{\cdots},\mu_{r_t},\psi\phi(\mu_{r_t})\}$ (cancel $\psi\phi(\mu_{r_i})$ if $\psi\phi(\mu_{r_i})=\mu_{r_i}$). So $|\Psi||\Phi|$ is also homotopic to the identity map of $|S(2^{\o K\setminus\{\phi\}})|$ by linear homotopy.

Let  $K_{\sigma,\omega}$ be as in Definition 1.2 and $\frak S_{\sigma,\omega}$ be as defined in Definition 3.4. Define maps $\psi_{\sigma,\omega}\colon K_{\sigma,\omega}{\setminus}\{\phi\}\to \frak S_{\sigma,\omega}{\setminus}\{\phi\}$ and $\varphi_{\sigma,\omega}\colon \frak S_{\sigma,\omega}{\setminus}\{\phi\}\to K_{\sigma,\omega}{\setminus}\{\phi\}$ as follows. For a non-empty simplex $\eta$ of $K_{\sigma,\omega}$, $\psi_{\sigma,\omega}(\eta)=\{k\,|\,\sigma{\cup}\eta\subset \sigma_{k}\}$. For a non-empty simplex $\mu=\{i_0,{\cdots},i_s\}$ of $\frak S_{\sigma,\omega}$, $\varphi_{\sigma,\omega}(\mu)=({\cap}_{j=0}^s\sigma_{i_j}){\setminus}\sigma$. $\psi_{\sigma,\omega}$ and $\varphi_{\sigma,\omega}$ are order reversing and so are not simplicial maps. But they induce simplicial maps $\Psi_{\sigma,\omega}\colon S(K_{\sigma,\omega})\to S(\frak S_{\sigma,\omega})$ and $\Phi_{\sigma,\omega}\colon S(\frak S_{\sigma,\omega})\to S(K_{\sigma,\omega})$ defined by $\Psi_{\sigma,\omega}(\{\eta_1,{\cdots},\eta_s\})
=\{\psi_{\sigma,\omega}(\eta_s),{\cdots},\psi_{\sigma,\omega}(\eta_1)\}$ and $\Phi_{\sigma,\omega}(\{\mu_1,{\cdots},\mu_t\})
=\{\varphi_{\sigma,\omega}(\mu_t),{\cdots},\varphi_{\sigma,\omega}(\mu_1)\}$, where repetition is canceled in the image set. $S(\frak S_{\sigma,\omega})$ is naturally a simplicial subcomplex of $S(2^{\o K\setminus\{\phi\}})$. Regard $S(K_{\sigma,\omega})$ as a simplicial subcomplex of $S(2^{K\setminus\{\phi\}})$ by the correspondence $\{\eta_1,\cdots,\eta_s\}\to\{\sigma{\cup}\eta_1,\cdots,\sigma{\cup}\eta_s\}$. Denote both the inclusion map by $i$. Then by definition, $i\,\Psi_{\sigma,\omega}=\Psi\, i$ and $i\,\Phi_{\sigma,\omega}=\Phi\, i$. This implies that $V$ and $\Phi_{\sigma,\omega}\Psi_{\sigma,\omega}(V)$ are faces of the same simplex for all non-empty simplex $V$ of $S(K_{\sigma,\omega})$ and $|\Phi_{\sigma,\omega}||\Psi_{\sigma,\omega}|$ is homotopic to the identity map of $S(K_{\sigma,\omega})$. Similarly, $|\Psi_{\sigma,\omega}||\Phi_{\sigma,\omega}|$ is homotopic to the identity map of $S(\frak S_{\sigma,\omega})$.
So we have the following commutative diagram of long exact sequences of simplicial homology groups
\[\begin{array}{ccccccccc}
\cdots\to&H_k(S(\frak S_{\sigma,\omega}))&\stackrel{i}{\to}&H_k(S(2^{\o K\setminus\{\phi\}}))&\to&
H_k(S(2^{\o K\setminus\{\phi\}}),S(\frak S_{\sigma,\omega}))&\to&H_{k-1}(S(K_{\sigma,\omega}))&\to\cdots\,\\
&{\scriptstyle \Phi_{\sigma,\omega}}\downarrow\quad&&{\scriptstyle \Phi}\downarrow\quad&&{\scriptstyle \o\Phi_{\sigma,\omega}}\downarrow
&&{\scriptstyle \Psi_{\sigma,\omega}}\downarrow\quad&\vspace{1mm}\\
{\cdots}\to & H_k(S(K_{\sigma,\omega})) &\stackrel{i}{\to}& H_k(S(2^{K\setminus\{\phi\}})) & \to &
H_k(S(2^{K\setminus\{\phi\}}),S(K_{\sigma,\omega}))& \to & H_{k-1}(S(K_{\sigma,\omega}))& \to {\cdots},
\end{array}\]
where we use the same symbol to denote a simplicial map and the homology homomorphism induced by it. Since both $\Phi$ and $\Phi_{\sigma,\omega}$ are isomorphisms, $\o\Phi_{\sigma,\omega}$ is also an isomorphism.

Suppose $x_i=\mu_i{\otimes}\xx_i$ for $x_1,x_2,x_3$ in the theorem. For ordered subset $\{i_0,\cdots,i_s\}\subset\mu_1$, define \[\w\Delta_{x_1}^{x_2,x_3}(\{i_0,\cdots,i_s\}{\otimes}\xx_1)=\Sigma_{k=0}^s
\{i_0,\cdots,i_k\}{\otimes}\xx_2\otimes\{i_k,\cdots,i_s\}{\otimes}\xx_3,\]
where $\{j_1,{\cdots},j_u\}{\otimes}(x_1,{\cdots},x_m)=0$ if there is $l\in\cap_{v=1}^u\sigma_{j_v}$ such that $x_l\in\frak k_l$.
By the isomorphism in the proof of Theorem 2.5, we have the following commutative diagram
\[\begin{array}{ccc}
(C_{*,*}(x_1),d)
&\stackrel{\w\Delta_{x_1}^{x_2,x_3}}{-\!\!\!-\!\!\!-\!\!\!-\!\!\!\longrightarrow}
&(C_{*,*}(x_2){\otimes}C_{*,*}(x_3),d{\otimes}d)\vspace{1mm}\\
\|\wr&&\|\wr\\
(C_*(2^{\mu_1},\frak S_{\sigma_1,\omega_1}),d)
&\stackrel{\w\Delta_{1}^{2,3}}{-\!\!\!-\!\!\!-\!\!\!-\!\!\!\longrightarrow}
&(C_*(2^{\mu_2},\frak S_{\sigma_2,\omega_2}){\otimes}C_*(2^{\mu_3},\frak S_{\sigma_3,\omega_3}),d{\otimes}d)  .\end{array}
\]

We have chain complex homomorphism $\Delta\colon (C_*(2^{[n]}),d)\to(C_*(2^{[n]}){\otimes}C_*(2^{[n]}),d{\otimes}d)$ defined as follows. For $\{i_0,\cdots,i_s\}\in C_*(2^{[n]})$, $\Delta(\{i_0,\cdots,i_s\})
=\Sigma_{k=0}^s\{i_0,\cdots,i_k\}{\otimes}\{i_k,\cdots,i_s\}$. If $(\cap_{j=0}^s\sigma_{i_j})\cap\omega_1\neq\phi$, then at least one of $(\cap_{j=0}^k\sigma_{i_j})\cap\omega_2$ and $(\cap_{j=k}^s\sigma_{i_j})\cap\omega_3$ is non-empty. This implies that\\
\hspace*{32mm}$\Delta(C_*(\frak S_{\sigma_1,\omega_1}))\subset C_*(\frak S_{\sigma_2,\omega_2}){\otimes}C_*(2^{[n]}){+}C_*(2^{[n]}){\otimes}C_*(\frak S_{\sigma_3,\omega_3})\vspace{1mm}$.\\
$\sigma_2\subset\sigma_1$ and $\sigma_3\subset\sigma_1$ imply that
$C_*(2^{\mu_1})\subset C_*(2^{\mu_2})$ and $C_*(2^{\mu_1})\subset C_*(2^{\mu_3})$.

So we have the following two commutative diagrams of short exact sequences of simplicial chain complexes ($C_i=C_*(\frak S_{\sigma_i,\omega_i})$, $\w C_i=C_*(2^{\mu_i})$ and $C=C_*(2^{[n]})$)
\[\begin{array}{ccccccc}
0\to&C_1&\to&C&\to&C/C_1&\to 0 \\
&{\scriptstyle\Delta}\downarrow\quad&&{\scriptstyle\Delta}\downarrow\quad
&&{\scriptstyle\Delta_1^{2,3}}\downarrow\quad &\vspace{1mm}\\
0\to&C_2{\otimes}C{+}C{\otimes}C_3&
\to&C{\otimes}C&\to&(C_*/C_2){\otimes}
(C/C_3)&\to 0  \end{array}
\]
\[\begin{array}{ccccccc}
0\to&C_1&\to&\w C_1&\to&\w C_1/C_1&\to 0 \\
&{\scriptstyle\Delta}\downarrow\quad&&{\scriptstyle\Delta}\downarrow\quad
&&{\scriptstyle\w\Delta_{1}^{2,3}}\downarrow\quad &\vspace{1mm}\\
0\to&C_2{\otimes}\w C_3{+}\w C_2{\otimes}C_3&
\to&\w C_2{\otimes}\w C_3&\to&(\w C_2/C_2){\otimes}
(\w C_3/C_3)&\to 0  \end{array}
\]
and from the above diagrams, we have the following commutative diagram
\[\begin{array}{ccc}
C_*(2^{\mu_1},\frak S_{\sigma_1,\omega_1})&\stackrel{\w\Delta_{1}^{2,3}}{-\!\!\!-\!\!\!-\!\!\!\longrightarrow}&
C_*(2^{\mu_2},\frak S_{\sigma_2,\omega_2}){\otimes}
C_*(2^{\mu_3},\frak S_{\sigma_3,\omega_3})\vspace{1mm}\\
\|\wr\quad&&\|\wr\quad\\
C_*(2^{[n]},\frak S_{\sigma_1,\omega_1})&\stackrel{\Delta_1^{2,3}}{-\!\!\!-\!\!\!-\!\!\!\longrightarrow}&
C_*(2^{[n]},\frak S_{\sigma_2,\omega_2}){\otimes}
C_*(2^{[n]},\frak S_{\sigma_3,\omega_3})\vspace{1mm}\\
\|\wr\quad&&\|\wr\quad\\
C_*(2^{\o K\setminus\{\phi\}},\frak S_{\sigma_1,\omega_1})&\stackrel{\Delta_1^{2,3}}{-\!\!\!-\!\!\!-\!\!\!\longrightarrow}&
C_*(2^{\o K\setminus\{\phi\}},\frak S_{\sigma_2,\omega_2}){\otimes}
C_*(2^{\o K\setminus\{\phi\}},\frak S_{\sigma_3,\omega_3})\vspace{1mm}\\
\|\wr\quad&&\|\wr\quad\\
C_*(S(2^{\o K\setminus\{\phi\}}),S(\frak S_{\sigma_1,\omega_1}))&\stackrel{\Delta_1^{2,3}}{-\!\!\!-\!\!\!-\!\!\!\longrightarrow}&
C_*(S(2^{\o K\setminus\{\phi\}}),S(\frak S_{\sigma_2,\omega_2})){\otimes}
C_*(S(2^{\o K\setminus\{\phi\}}),S(\frak S_{\sigma_3,\omega_3}))\end{array}
\]

Apply the isomorphism $\o\Phi_{\sigma_i,\omega_i}\colon C_*(S(2^{\o K\setminus\{\phi\}}),S(\frak S_{\sigma_i,\omega_i}))\to C_*(S(2^{K\setminus\{\phi\}}),S(K_{\sigma_i,\omega_i}))$ to the last homomorphism, we get a chain complex homomorphism
\[C_*(S(2^{K\setminus\{\phi\}}),S(K_{\sigma_1,\omega_1}))
\stackrel{\Delta_1^{2,3}}{-\!\!\!-\!\!\!-\!\!\!\longrightarrow}
C_*(S(2^{K\setminus\{\phi\}}),S(K_{\sigma_2,\omega_2})){\otimes}
C_*(S(2^{K\setminus\{\phi\}}),S(K_{\sigma_3,\omega_3})).\]

Notice that all $\Delta_1^{2,3}$ is induced by a simplicial approximation of the diagonal map $\Delta$ defined by $\Delta(x)=(x,x)$. Since barycentric subdivision does not change the homology, we have the following commutative diagram
\[\begin{array}{ccc}
H_*(S(2^{K\setminus\{\phi\}}),S(K_{\sigma_1,\omega_1}))
&\stackrel{(\Delta_1^{2,3})_*}{-\!\!\!-\!\!\!-\!\!\!\longrightarrow}&
H_*(S(2^{K\setminus\{\phi\}}),S(K_{\sigma_2,\omega_2})){\otimes}
H_*(S(2^{K\setminus\{\phi\}}),S(K_{\sigma_3,\omega_3}))\\
\|\wr\quad&&\|\wr\quad\\
H_*(|S(2^{K\setminus\{\phi\}})|,|S(K_{\sigma_1,\omega_1})|)&
\stackrel{S(\Delta)_*}{-\!\!\!-\!\!\!-\!\!\!\longrightarrow}&
H_*(|S(2^{K\setminus\{\phi\}})|,|S(K_{\sigma_2,\omega_2})|){\otimes}
H_*(|S(2^{K\setminus\{\phi\}})|,|S(K_{\sigma_3,\omega_3})|)\vspace{1mm}\\
\|\wr\quad&&\|\wr\quad\\
H_*(|2^{K\setminus\{\phi\}}|,|K_{\sigma_1,\omega_1}|)&
\stackrel{\Delta_*}{-\!\!\!-\!\!\!-\!\!\!\longrightarrow}&
H_*(|2^{K\setminus\{\phi\}}|,|K_{\sigma_2,\omega_2}|){\otimes}
H_*(|2^{K\setminus\{\phi\}}|,|K_{\sigma_3,\omega_3}|)\vspace{1mm}\\
\|\wr\quad&&\|\wr\quad\vspace{1mm}\\
H_*(|2^{v_1}|,|K_{\sigma_1,\omega_1}|)&
\stackrel{\Delta_{v_1}^{v_2,v_3}}{-\!\!\!-\!\!\!-\!\!\!\longrightarrow}&
H_*(|2^{v_2}|,|K_{\sigma_2,\omega_2}|){\otimes}
H_*(|2^{v_3}|,|K_{\sigma_3,\omega_3}|)\end{array}\]
So $\Delta_{x_1}^{x_2,x_3}$ and $\Delta_{v_1}^{v_2,v_3}$ satisfy the commutative diagram of the theorem.\end{proof}\vspace{3mm}

{\bf Definition 3.7} {\it Let $M=\ZZ_K(X,A)$ and $T_{*}(M)$ be as in Definition 3.2. We use the same symbol to denote both the base and the space generated by the base. Then the coproduct $\Delta^T\colon T_{*}(M)\to T_{*}(M){\otimes}T_{*}(M)$ is defined as follows. For $\xx=(x_1,\cdots,x_m)\in T_{*}(M)$, regard $x_k$ as an element of $H_*(A_k)$ if $x_k\in\frak i_k$. Suppose $\Delta_k(x_k)=\Sigma_jc_{k,j}x'_{k,j}{\otimes}x''_{k,j}$, where $\Delta_k$ is the coproduct of the coalgebra $H_*(X_k)$ or $H_*(A_k)$, $c_{k,j}\in\kk{\setminus}\{0\}$, and $x'_{k,j},x''_{k,j}$ are base elements of $H_*(X_k)$ or $H_*(A_k)$. Then
\[\Delta^T(\xx)=\Sigma_{j_1,\cdots,j_m}(-1)^\varepsilon c_{1,j_1}\cdots c_{m,j_m} (x'_{1,j_1},\cdots,x'_{m,j_m}){\otimes}(x''_{1,j_1},\cdots,x''_{m,j_m}),\]
where $\varepsilon=|x''_{1,j_1}|(|x'_{2,j_2}|{+}{\cdots}{+}|x'_{m,j_m}|){+}
|x''_{2,j_2}|(|x'_{3,j_3}|{+}{\cdots}{+}|x'_{m,j_m}|){+}{\cdots}{+}|x''_{m-1,j_{m-1}}||x'_{m,j_m}|$.
}\vspace{3mm}

Notice that in general, $\Delta^T$ is not coassociative.\vspace{3mm}

{\bf Theorem 3.8} {\it Let everything be as in Theorem 1.3 and Definition 3.7. The coproduct of \[H_*(M)=\oplus_{(\sigma,\omega)\in I_M} H_*^{\sigma,\omega}(M)=\oplus_{(\sigma,\omega)\in I_M}H_*(|2^{v}|,|K_{\sigma,\omega}|){\otimes}T_*^{\sigma,\omega}(M)\] is as follows. For $\xx\in T_*^{\sigma,\omega}(M)$ and $a\in H_*(|2^v|,|K_{\sigma,\omega}|)$, if $\Delta^T(\xx)=\Sigma_jc_j\xx'_j{\otimes}\xx''_j$ and $\Delta_v^{v'_j,v''_j}(a)=\Sigma_{k}a'_{j,k}{\otimes}a''_{j,k}$, where $\xx'_j\in T_*^{\sigma'_j,\omega'_j}(M)$ and  $\xx''_j\in T_*^{\sigma''_j,\omega''_j}(M)$, $v'_j$ and $v''_j$ are respcetively the vertex set of $K_{\sigma'_j,\omega'_j}$ and $K_{\sigma''_j,\omega''_j}$, $\Delta_v^{v'_j,v''_j}$ is as defined in Theorem 3.6, then\vspace{1mm}\\
\hspace*{43.5mm}$\Delta(a{\otimes}\xx)=\Sigma_{j,k}\,(-1)^{|x'_j||a''_{j,k}|}c_j(a'_{j,k}{\otimes}\xx'_j)\otimes (a''_{j,k}{\otimes}\xx''_j).$

Dually, the product of \[H^*(M)=\oplus_{(\sigma,\omega)\in I_M} H^*_{\sigma,\omega}(M)=\oplus_{(\sigma,\omega)\in I_M}H^*(|2^{v}|,|K_{\sigma,\omega}|){\otimes}T^*_{\sigma,\omega}(M)\] is as follows. For $\xx_i\in T^*_{\sigma_i,\omega_i}(M)$ and $a_i\in H^*(|2^{v_i}|,|K_{\sigma_i,\omega_i}|)$, $i=1,2$, if $\langle\xx_1,\xx_2\rangle_T=\Sigma_jc_j\xx_j$ with $\xx_j\in T^*_{\sigma_j,\omega_j}(M)$ and $v_j$ the vertex set of $K_{\sigma_j,\omega_j}$, where $\langle\,,\rangle_T$ is the dual product of $\Delta^T$ in Definition 3.7, then\vspace{1mm}\\
\hspace*{38mm}$(a_1{\otimes}\xx_1)(a_2{\otimes}\xx_2)=\Sigma_j(-1)^{|\xx_1||a_2|}c_j\langle a_1,a_2\rangle_{v_1,v_2}^{v_j}{\otimes}\xx_j$,\\
where $\langle\,,\rangle_{v_1,v_2}^{v_j}$ is the dual product of $\Delta_{v_1}^{v_2,v_3}$ in Theorem 3.6.
}
\begin{proof} Suppose $K=\{\sigma_1,\cdots,\sigma_n\}$ (in any order), then $\frak S=(\sigma_1,\cdots,\sigma_n)$ is a simplex cover of $K$. Let $\frak C$ be the cover associated to $\frak S$. Consider the coproduct $\Delta\colon\frak C_{*,*}(M)\to\frak C_{*,*}(M){\otimes}\frak C_{*,*}(M)$ in Definition 2.2. Let everything be as in Theorem 3.3 and its proof. For $x=\mu{\otimes}\xx=\mu{\otimes}(x_1,\cdots,x_m)\in T_{*,*}^{\sigma,\omega}(\frak C)$, let $\eta=\cap_{j\in\mu}\sigma_j$. Now we prove that there is no $k\in\eta$ such that $x_k\in\frak i_k$. Suppose there is $k\in\eta$ such that $x_k\in\frak i_k$. Let $\eta'=\eta{\setminus}\{k\}$, then $\eta'\in K$ and there is a $j\in[n]$ such that $\eta'=\sigma_j$. $\mu{\cup}\{j\}{\otimes}\xx$ is in the same connected component of $x$. A contradiction to that $x$ is a top vertex! So there is no $k\in\eta$ such that $x_k\in\frak i_k$. This implies that $\Delta_\mu(\xx)=\Delta^T(\xx)$, where $\Delta_\mu$ is the coproduct of $H_*(M_\mu)$ as in Definition 3.2. Suppose $\Delta^T(\xx)=\Sigma_jc_j\xx'_j{\otimes}\xx''_j$, where $c_j\in\kk{\setminus}\{0\}$ and $\xx'_j,\xx''_j\in T_*(M)$. Then for ordered subset $\{i_0,\cdots,i_s\}\subset\mu$,
\begin{eqnarray*}&&\Delta(\{i_0,\cdots,i_s\}{\otimes}\xx)\\
&=&\Sigma_j\Sigma_{k=0}^s(-1)^{(s-k)|\xx'_j|}c_j\big(\{i_0,\cdots,i_k\}
{\otimes}\xx'_j\big)\otimes\big(\{i_k,\cdots,i_s\}{\otimes}\xx''_j\big),
\end{eqnarray*}
where $\{j_1,{\cdots},j_u\}{\otimes}(x_1,{\cdots},x_m)=0$ if there is $l\in\cap_{v=1}^u\sigma_{j_v}$ such that $x_l\in\frak k_l$.
Compare $\Delta(\{i_0,\cdots,i_s\}{\otimes}\xx)$ with $\w\Delta_{x_1}^{x_2,x_3}(\{i_0,\cdots,i_s\}{\otimes}\xx_1)$ defined in the proof of Theorem 3.6. They differ only in sign. So the induced homology homomorphism (still denoted by) $\Delta$ is just the formula in the theorem.
\end{proof}\vspace{3mm}

{\bf Theorem 3.9} {\it Let everything be as in Theorem 3.8. If for every $k$, there is a subcoalgebra $\frak i_k$ of $H_*(A_k)$ such that the restriction of $i_k$ on $\frak i_k$ is a coalgebra isomorphism from ${\rm coim}i_k$ to ${\rm im}i_k$, then $(T_*(M),\Delta^T)$ is a cocommutative, coassociative (in the graded sense) coalgebra and $(T^*(M),\langle\,,\rangle_T)$ is a commutative, associative algebra over $\kk$.}
\begin{proof} In general, ${\rm coim}i_k={\rm im}i_k$ only as Abellian groups, the coproduct of the two groups are different. But in the condition of the theorem, ${\rm coim}i_k={\rm im}i_k$ as coalgebras. So the cocommutativity and coassociativity of $\Delta^T$ follow from that of every $H_*(A_k)$ and $H_*(X_k)$.
\end{proof}

{\bf Remark} When $M=\ZZ_K(D^2,S^1)$, the usual moment-angle complex with every space pair $(X_k,A_k)=(D^2,S^1)$, there is an algebra isomorphism $H^*(M)={\rm Tor}^{\kk[\xx]}_*(F(K),\kk)$, where $F(K)$ is the Stanley-Reisner face ring of $K$. But the formula even can not be generalized to the following most similar case when every $\w H_*(X_k)=0$ and $\w H_s(A_k)=\kk$ for $s=r_k>1$ and $\w H_s(A_k)=0$ otherwise. Since the degree $r_k$ may be even or odd and the degree of $x_k$ can only be even, the formula can not be generalized to this case by giving the polynomial ring $\kk[x_1,\cdots,x_m]$ a graded commutative algebra structure.\vspace{3mm}

We end this paper with an example. Let $M=\ZZ_K\Big(\begin{array}{ccc}
\scriptstyle{r_1{+}1} &\scriptstyle{\cdots}&\scriptstyle{r_m{+}1}\\
\scriptstyle{ k_1}   &\scriptstyle{\cdots}& \scriptstyle{k_m}\end{array}\Big)$ be as in Example 1.7. It is a generalized moment-angle complex that satisfies the condition of Theorem 3.9. Since all $T_*^{\sigma,\omega}(M)$ is one dimensional, we have
\[H_*^{\sigma,\omega}(M)=H_*(|2^{v}|,|K_{\sigma,\omega}|){\otimes}(n_1,\cdots,n_m)
=\Sigma^{|(\sigma,\omega)|}H_*(|2^{v}|,|K_{\sigma,\omega}|),\]
where $\Sigma^r$ means uplift the degree by $r$, $|(\sigma,\omega)|=\Sigma_{i\in\sigma}(r_i{+}1)+\Sigma_{i\in\omega}k_i$,  $n_i=r_i{+}1$ if $i\in\sigma$, $n_i=0$ if $i\in\w\sigma$, $n_i=k_i$ if $i\in\omega$. $(\sigma,\omega)\to (n_1,\cdots,n_m)$ is a 1-1 correspondence from $I_M$ to $T_*(M)$. Identify the two sets and the coproduct $\Delta^T$ is defined by that for all $(\sigma,\omega)\in I_M$\\
\hspace*{55mm}$\Delta^T(\sigma,\omega)=\Sigma\, (-1)^\varepsilon(\sigma',\omega'){\otimes}(\sigma'',\omega'')$\\
where $\varepsilon=n''_1(n'_2{+}{\cdots}{+}n'_m)+
n''_2(n'_3{+}{\cdots}{+}n'_m)+{\cdots}+n''_{m-1}n'_m$ and the sum is taken over all $\sigma'{\sqcup}\sigma''=\sigma$, $\omega'{\sqcup}\omega''=\omega$ ($\sqcup$ means disjoint union). Suppose $v',v''$ are respectively the vertex set of $K_{\sigma',\omega'}$ and $K_{\sigma'',\omega''}$ and $\Delta_v^{v',v''}\colon H_*(|2^{v}|,|K_{\sigma,\omega}|)\to
H_*(|2^{v'}|,|K_{\sigma',\omega'}|){\otimes}H_*(|2^{v''}|,|K_{\sigma'',\omega''}|)$ is as defined in Theorem 3.6. Then for $a\in H_*(|2^{v}|,|K_{\sigma,\omega}|)$ with $\Delta_v^{v',v''}(a)=\Sigma_j a'_j{\otimes}a''_j$,
\[\Delta(a)=\sum_{\sigma'\sqcup\sigma''=\sigma,\, \omega'\sqcup\omega''=\omega}\sum_j\,(-1)^{\varepsilon+|a''_j||(\sigma',\omega')|}
a'_j{\otimes}a''_j,\]
where the element $b{\otimes}(l_1,\cdots,l_m)$ in $H_*(|2^{w}|,|K_{\rho,\varrho}|){\otimes}(l_1,\cdots,l_m)$ is still simply denoted by $b$.

Dually,
\[H^*(M)=\oplus_{(\sigma,\omega)\in I_M}\Sigma^{|(\sigma,\omega)|}H^*(|2^{v}|,|K_{\sigma,\omega}|).\]
For $a'\in H^*(|2^{v'}|,|K_{\sigma',\omega'}|)$ and
$a''\in H^*(|2^{v''}|,|K_{\sigma'',\omega''}|)$, still denote the corresponding element in $H^*_{\sigma',\omega'}$ and $H^*_{\sigma'',\omega''}$ by $a'$ and $a''$. Then $a'a''=0$ except that $\sigma'{\cap}\sigma''=\phi$, $\omega'{\cap}\omega''=\phi$ and $(\sigma'{\cup}\sigma'',\omega'{\cup}\omega'')\in I_M$. For such $\sigma',\sigma''$ and $\omega',\omega''$,
\[a'a''=\sum\,(-1)^{\varepsilon+|a''||(\sigma',\omega')|}
\langle a',a''\rangle_{v',v''}^v\]
where $\varepsilon$ is as above and $\langle\,,\rangle^v_{v',v''}\colon
H^*(|2^{v'}|,|K_{\sigma',\omega'}|){\otimes}H^*(|2^{v''}|,|K_{\sigma'',\omega''}|)\to
H^*(|2^{v}|,|K_{\sigma'\sqcup\sigma''\!,\omega'\sqcup\omega''}|)$ is the dual map of $\Delta_v^{v',v''}$.


\vspace{6mm}


\begin{thebibliography}{00}

\bibitem{AP} C.~Allday and V.~Puppe, {\it Cohomological Methods in Transformation Groups,}
 Cambridge Studies in Advanced Mathematics, {\bf 32}, Cambridge University Press, 1992.
\bibitem{BB} A.~Bahri, M.Bendersky, F.~R.~Cohen and S.~Gitler, Decompositions of the polyhedral
 product functor with applications to moment-angle complexes and related spaces,
 {\it Proc. Nat. Acad. Sci. U. S. A.} {\bf 106} (2009), 12241-12243.
\bibitem{BBCG08} A.~Bahri, M.~Bendersky, F.~R.~Cohen and S.~Gitler,
{\it The polyhedral product functor: A methord of computation for moment-angle complexes,
arrangement and related spaces,}  preprint arXiv:0711.4689v2 [math.AT] 8 Dec 2008.
\bibitem{Baskakov} I.~Baskakov, {\it Cohomology of K-powers of spaces and the combinatorics
of simplicial divisions,} Russian Math. Surveys {\bf 57} (2002), no. 5, 989-990.
\bibitem{BP} V.~M.~Buchstaber and T.~E.~Panov, {\it Torus Actions and Their Applications in
 Topology and Combinatorics}, University Lecture Series, Vol. {\bf 24}, Amer. Math. Soc.
 Providence, RI, 2002.
\bibitem{BP04} V.~M.~Buchstaber and T.~E.~Panov, Combinatorics of simplicial cell complexes
 and tours action, {\it Proc. Steklov Inst. Math.} {\bf 247} (2004), 33-49.
\bibitem{CL} X.~Cao and Z.~L\"{u},  M\"{o}bius transform, moment-angle complexes and
 Halperin-Carlsson conjecture, {\it preprint.} arXiv:0908.3174v2 [math.CO] 12 Sep 2009.
\bibitem{DJ} M.~W.~Davis and T.~Januszkiewicz,  Convex polytopes, coxeter orbifolds and
 torus action, {\it Duck Math. J.} {\bf 62} (1991), 417-451.
\bibitem{DS} G.~Denham and A.~Suciu, Moment-angle complexes, monomial ideals and Massey products,
 {\it Pure Appl. Math. Q.} {\bf 3} (2007), 25-60.
\bibitem{Franz} M.~Franz, {The intergral cohomology of toric manifolds,} Proc. Steklov Inst. Math.
{\bf 252} (2006), 53-62. [Proceedings of the Keldysh Conference, Moscow 2004].
\bibitem{GM} R.~Goresky and R.~MacPherson, {\it Stratified Morse Theory,} Ergeb. Math. Grenzgeb.,
Vol. 14, Springer-Verlag, Berlin, 1988.
\bibitem{GTh} J.~Grbic and R.~Theriault, {\it Homotopy type of the complement of a coordinate
sunspace arrangement of codimension two,} Russian Math. Surveys {\bf 59} (2004), no. 3, 1207-1209.
\bibitem{Hilton} P.~J.~Hilton, U.~Stammbach, {\it A Course in Homological Algebra,}
 Berlin-Heideberg-New York: Springer 1971.
\bibitem{Hoc} M.~Hochster, {\it Cohen-Macaulay ring, combinatorics and simplecial complexes,}
in: Ring theory, II (Proc. Second Conf. Univ. Oklahoma, Norman, Okla., 1975), pp. 171-222.
\bibitem{Lopez} S.~Lopez de Medrano, {\it Topology of the intersection of quadrics in $\mathbb{R}^n$,}
in {\it Algebraic Topology} (Arcata Ca), Springer Verlag LNM {\bf 1370} (1989), Springer Verlag.
\bibitem{LP} Z.~Lu and T.~Panov, Moment-angle complexes from simplicial posets, {\it preprint}
arXiv:0912.2219v1 [math.AT] 11 Dec 2009.
\bibitem{Mil} E.~Miller and B.~Sturmfels, {\it Combinatorial Commutative Algebra},
 Graduate Texts in Math. {\bf 227}, Springer, 2005.
\bibitem{Panov} T.~E.~ Panov, Cohomology of face rings and tours actions, {\it London
 Math. Soc. Lect. Notes Ser.} {\bf 347} (2008), 165-201.
\bibitem{Porter} G.~Porter, {\it The homotopy groups of wedge of suspensions,} Amer. J. Math.,
{\bf 88} (1966), 655-662.
\bibitem{Stanley} R.~P.~Stanley {\it Combinatorics and Commutative Algebra} second edition, Progress in Math.
{\bf 41} Birkhauser, Boston, 1996.
\bibitem{U} Y.~Ustinovsky, Toral rank conjecture for moment-angle complexes, {\it preprint}
arXiv:0909.1053v2 [math.AT] 29 Sep 2009.
\bibitem{Vogt} R.~Vogt, {\it Homotopy limits and colimits,} Math. Z. {\bf 134} (1973), 11-52.
\bibitem{Zha} Q.Zheng, {\it A new Massey product on Ext groups},  Journal of Algebra., {\bf 183} (1996), 378-395.
\bibitem{Zhb} Q.Zheng, {\it S-module and the new Massey-product},  Journal of Algebra., {\bf 190} (1997), 478-397.
\bibitem{Z} Q. Zheng and X. Wang, The homology of simplicial complexes and the cohomology of moment-angle complexes {\it preprint}
arXiv:1109.6382v1[math.AT]
\end{thebibliography}
\end{document}